\documentclass[10pt,oneside,a4paper]{article}
\usepackage{times}
\usepackage{amscd}
\usepackage[latin1]{inputenc}
\usepackage{latexsym}         

\usepackage{amsmath,amssymb,amsthm,amsfonts,mathrsfs}

\usepackage{calc}
\usepackage{psfrag}
\usepackage{xypic}

\usepackage{graphicx,color}

\usepackage{epsf}
\newtheorem{teor}{Theorem}[section]
\newtheorem{lemma}[teor]{Lemma}
\newtheorem{propos}[teor]{Proposition}
\newtheorem{corol}[teor]{Corollary}
\newtheorem{claim}[teor]{Claim}
\theoremstyle{definition}
\newtheorem{defin}[teor]{Definition}
\newtheorem{osserv}[teor]{Remark}
\newtheorem{example}[teor]{Example}
\newcommand{\B}{\mbox{ $ \mathbb{B}$}}
\newcommand{\R}{\mbox{ $ \mathbb{R}$}}
\newcommand{\Q}{\mbox{ $ \mathbb{ Q}$}}
\newcommand{\dex}{\mbox{ $ \frac{\partial}{\partial x} $ }}
\newcommand{\dey}{\mbox{ $ \frac{\partial}{\partial y} $ }}
\newcommand{\h}{\mbox{ $  \mathbb {H} $ }}
\newcommand{\C}{\mbox{ $ \mathcal{C}$}}
\newcommand{\Pal}{\mbox{ $ \mathbb{P} $ }}
\newcommand{\Cool}{\mbox{ $ \mathbb{C}$}}
\newcommand{\N}{\mbox{ $\mathbb{N} $ }}
\newcommand{\Z}{\mbox{ $\mathbb{Z} $ }}
\newcommand{\Fl}{\mbox{$\mathscr{F}$}}
\newcommand{\Ba}{\mbox{$\mathcal{B}$}}
\newcommand{\Bi}{\mbox{$\mathscr{B}$}}
\newcommand{\Ol}{\mbox{ $ \mathcal {O}$}}
\newcommand{\Sal}{\mbox{ $ \mathcal{S} $}}
\newcommand{\Il}{\mbox{ $ \mathcal {I}$}}
\newcommand{\M}{\mbox{ $ M$ }}
\newcommand{\W}{\mbox{$\mathcal{W}$}}
\newcommand{\s}{\mbox{$ S$ }}
\newcommand{\p}{\mbox{$p$}}
\newcommand{\na}{\mbox{ $\nabla$}}
\newcommand{\f}{\mbox{ $ f $}}
\newcommand{\fu}{\mbox{ $ f:M \longrightarrow M $ }}
\newcommand{\id}{\mbox{  \hbox{Ind}}}
\newcommand{\F}{\mbox{$ \mathcal{F} $}}
\newcommand{\spec}{{\em{of the theorem}}}
\newcommand{\Pl}{\mbox{ $\mathcal{P}$}}
\newcommand{\bo}{\mbox{\textrm{B}}}
\newcommand{\T}{\mbox{ $  \mathbb {T} $ }}
\theoremstyle{remark}  %\newtheorem{nota}[teo]{Remark}
\numberwithin{equation}{section}

\title{On smooth foliations with Morse singularities}
\author{Lilia Rosati}
\date{Universit\`a di Firenze,\\Dipartimento di Matematica ``U. Dini'',\\viale
Morgagni 67/A, 50134 Firenze\\e-mail: {\small{{\texttt{rosati@math.unifi.it}}}}}
\begin{document}
\maketitle
\section*{Abstract}
Let $M$ be a smooth manifold and let $\F$ be a codimension one,
$C^\infty$ foliation on $M$, with isolated singularities of Morse type. The study and
classification of pairs $(M,\F)$ is a challenging (and difficult) problem. In this setting, a classical result due to Reeb \cite{Reeb} states that a manifold admitting a foliation with exactly two center-type singularities is a sphere. In particular this is true if the foliation is given by a function. Along these lines a result due to Eells and Kuiper \cite{Ku-Ee} classify manifolds having a real-valued function admitting exactly three non-degenerate singular points. In the present paper, we prove a generalization of the above mentioned results. To do this, we first describe the possible arrangements of pairs of singularities and the corresponding codimension one invariant sets, and then we give an elimination procedure for suitable center-saddle and some saddle-saddle configurations (of consecutive indices).\\
In the second part, we investigate if other classical results, such as Haefliger and Novikov (Compact Leaf) theorems, proved for regular foliations, still hold true in presence of singularities. At this purpose, in the singular set, $Sing(\F)$ of the foliation $\F$, we consider {\em{weakly stable}} components, that we define as those components admitting a neighborhood where all leaves are compact.  If $Sing(\F)$ admits only weakly stable components, given by smoothly embedded curves diffeomorphic to $S^1$, we are able to extend Haefliger's theorem. Finally, the existence of a closed curve, transverse to the foliation, leads us to state a Novikov-type result.\\   
\section*{Acknoledgements} I am very grateful to prof. Bruno Sc\'ardua for proposing me such an interesting subject and for his valuable advice.  My hearthy good thanks to prof. Graziano Gentili for his suggestions on the writing of this article.
\section{Foliations and Morse Foliations}\label{inizio}    
{\bf{Definition \ref{inizio}.1}} A {\em{codimension $k$, foliated manifold}} $(M,\F)$ is a manifold $M$ with a differentiable structure, given by an atlas $\{(U_i, \phi_i)\}_{i \in I}$, satisfying the following properties:\\
(1) $\phi_i(U_i)= \bo ^{n-k} \times \bo ^k$;\\
(2) in $U_i \cap U_j \neq \varnothing$, we have $\phi_j \circ \phi_i^{-1}(x,y)=(f_{ij}(x,y),g_{ij}(y))$,\\
where $\{f_{ij}\}$ and $\{g_{ij}\}$ are families of, respectively, submersions and diffeomorphisms, defined on natural domains. Given a local chart ({\em{foliated chart}}) $(U, \phi)$, $\forall x \in \bo ^{n-k}$ and $y \in \bo ^k$, the set $\phi^{-1}(\cdot , y)$ is a {\em{plaque}} and the set $\phi^{-1}(x,\cdot)$ is a {\em{transverse section}}.

The existence of a foliated manifold $(M,\F)$ determines a partition of $M$ into subsets, the {\em{leaves}}, defined by means of an equivalence relation, each endowed of an intrinsic manifold structure. Let $x \in M$; we denote by $\F_x$ or $L_x$ the leaf of $\F$ through $x$. With the intrinsic manifold structure, $\F_x$ turns to be an immersed (but not embedded, in general) submanifold of $M$.\\
In an equivalent way, a foliated manifold $(M,\F)$ is a manifold $M$ with a collection of couples $\{(U_i,g_i)\}_{i \in I}$, where $\{U_i\}_{i \in I}$ is an open covering of $M$, $g_i:U_i \rightarrow \bo ^k$ is a submersion, $\forall i \in I$, and the $g_i$'s satisfy the cocycle relations, $g_i=g_{ij} \circ g_j$, $g_{ii}=id$, for suitable diffeomorphisms $g_{ij}:\bo ^k \rightarrow \bo^k$, defined when $U_i \cap U_j \neq \varnothing$. Each $U_i$ is said a {\em{foliation box}}, and $g_i$ a {\em{distinguished map}}. The functions $\gamma_{ij}=\textrm{d}g_{ij}$ are the transition maps \cite{Stee} of a bundle $N \F \subset TM$, normal to the foliation. More completely, there exists a G-structure on $M$ \cite{Law}, which is a reduction of the structure group $GL(n,\R)$ of the tangent bundle to the subgroup of the matrices $\left(\begin{array}{c|c}
A & B \\\hline 0 & C \end{array}\right)$, where $A \in GL(n-k,\R)$ and $C \in GL(k, \R)$.

A codimension one, $C^\infty$ foliation of a smooth manifold $M$, with isolated
singularities, is a pair 
$\F=(\F^*,Sing(\F))$, where $Sing(\F) \subset M$ is a discrete subset and
$\F^*$ is a codimension one, $C^\infty$ foliation (in the ordinary
sense) of $M^*=M \setminus Sing(\F)$. The {\em{leaves}} of $\F$
are the leaves of $\F^*$ and $Sing(\F)$ is the {\em{singular set}}
of $\F$. A point $p$ is a {\em{Morse singularity}} if there is a
$C^\infty$ function, $f_p:U_p \subset M \rightarrow \R$, defined
in a neighborhood $U_p$ of $p$, with a (single) non-degenerate critical
point at $p$ and such that $f_p$ is a local first integral of the foliation, i.e. the leaves of the restriction $\F|_{U_p}$ are the connected components of the level
 hypersurfaces of $f_p$ in $U_p \setminus \{p \}$. A Morse singularity $p$,
  of index $l$, is a {\em{saddle}}, if $0<l<n$ (where $n=\dim M$), and a {\em{center}},
  if $l=0,n$.
We say that the foliation $\F$ has a {\em{saddle-connection}} when there exists a
leaf accumulated by at least two distinct saddle-points.
A {\em{Morse foliation}} is a foliation with isolated singularities, whose singular
set consists of Morse singularities, and which has no saddle-connections. In this way if a Morse foliation has a (global) first integral, it is given by a Morse function.\\
Of course, the first basic example of a Morse foliation is indeed a
foliation defined by a Morse function on $M$. A less evident
example is given by the foliation depicted in figure \ref{less}.

In the literature, the orientability of a codimension $k$ (regular) foliation is determined by the orientability of the $(n-k)$-plane field tangent to the foliation, $x \rightarrow T_x \F_x$. Similarly transverse orientability is determined by the orientability of a complementary $k$-plane field. A singular, codimension one foliation, $\F$, is {\em{transversely orientable}} \cite{Ca-Sca} if it is given by the natural $(n-1)$-plane field associated to a one-form, $\omega \in \Lambda ^1(M)$, which is integrable in the sense of Frobenius. In this case, choosing a Riemannian metric on $M$, we may find a global vector field transverse to the foliation, $X=grad(\omega)$, $\omega X \geq 0$, and $\omega_x X_x=0$ if and only if $x$ is a singularity for the foliation ($\omega(x)=0$). A transversely orientable, singular foliation 
$\F$ of $M$ is a transversely orientable (regular) foliation $\F^*$ of $M^*$ in the sense of the classical definition. Viceversa, if $\F^*$ is transversely orientable, in general, $\F$ is not.

Thanks to the Morse Lemma \cite{Mil1}, Morse foliations reduce to few representative cases. On the other hand, Morse foliations describe a large class among transverseley orientable foliations. To see this, let $\F$ be a foliation defined by an integrable one-form, $\omega \in \Lambda^1(M)$, with isolated singularies. We proceed with a local analysis; using a local chart around each singularity, we may suppose $\omega \in \Lambda^1(\R^n)$, $\omega(0)=0$, and 0 is the only singularity of $\omega$. We have $\omega(x)= \sum_i h_i(x) dx^i$ and, in a neighborhood of $0 \in \R^n$, we may write $\omega(x)=\omega_1(x)+O(|x|^2)$, where $\omega_1$ is the linear part of $\omega$, defined by $\omega_1(x)= \sum_{i,j}a_{ij}x^i dx^j$, $a_{ij}=\partial h^i(x)/\partial x^j$. We recall that the integrability of $\omega$ implies the integrability of $\omega_1$ and that the singularity $0 \in \R^n$ is said to be non degenerate if and only if $(a_{ij}) \in \R(n)$ is non degenerate; in this latter case $(a_{ij})$ is symmetric: it is the hessian matrix of some real function $f$, defining the linearized foliation ($\omega_1= \textrm{d}f$). We have
\begin{displaymath}
\begin{array}{ccc}
\{\textrm{transverseley orientable foliations, with Morse singularities}\}=\\\{\textrm{foliations, defined by non degenerate linear one-forms}\} \subset \\\{\textrm{foliations, defined by non degenerate one-forms}\}.
\end{array}\end{displaymath}   
Let $(\sigma, \tau)$ be the space $\sigma$ of integrable one-forms in $\R^n$, with a singularity at the origin, endowed with the $C^1$-Whitney topology, $\tau$. If $\omega,\omega' \in \sigma$, we say $\omega$ {\em{equivalent}} $\omega'$ ($\omega \sim \omega'$) if there exists a diffeomorphism $\phi:\R^n \rightarrow \R^n$, $\phi(0)=0$, which sends leaves of $\omega$ into leaves of $\omega'$. Moreover, we say $\omega$ is {\em{structurally stable}}, if there exists a neighborhood $V$ of $\omega$ in $(\sigma,\tau)$ such that $\omega' \sim \omega, \forall \omega' \in V$.\\
{\bf{Theorem \ref{inizio}.2 (Wagneur)}}\cite{Wag} The one-form {\em{$\omega \in \sigma$ is structurally stable, if and only if the index of $0 \in Sing(\omega)$ is neither $2$ nor $n-2$}}.

Let us denote by $S$ the space of foliations defined by non degenerate one-forms with singularities, whose index is neither $2$ nor $n-2$. If $S_1 \subset S$ is the subset of foliations defined by linear one-forms, then we have:\\
{\bf{Corollary \ref{inizio}.3}} {\em{There exists a surjective map,}}$$s: S_1 \rightarrow S/_\sim.$$ 
\section{Holonomy and Reeb Stability Theorems}\label{uno}
\begin{figure}
\begin{minipage}[t]{.45\linewidth}
  {
\begin{center}
{
        \psfrag{1}{$\F_1$}
        \psfrag{2}{$\F_2$}
        \psfrag{l}{$L$}
        \psfrag{l0}{$L_0$}
        \psfrag{l1}{$L_1$}
        \psfrag{l2}{$L_2$}
        \includegraphics[scale=.4]{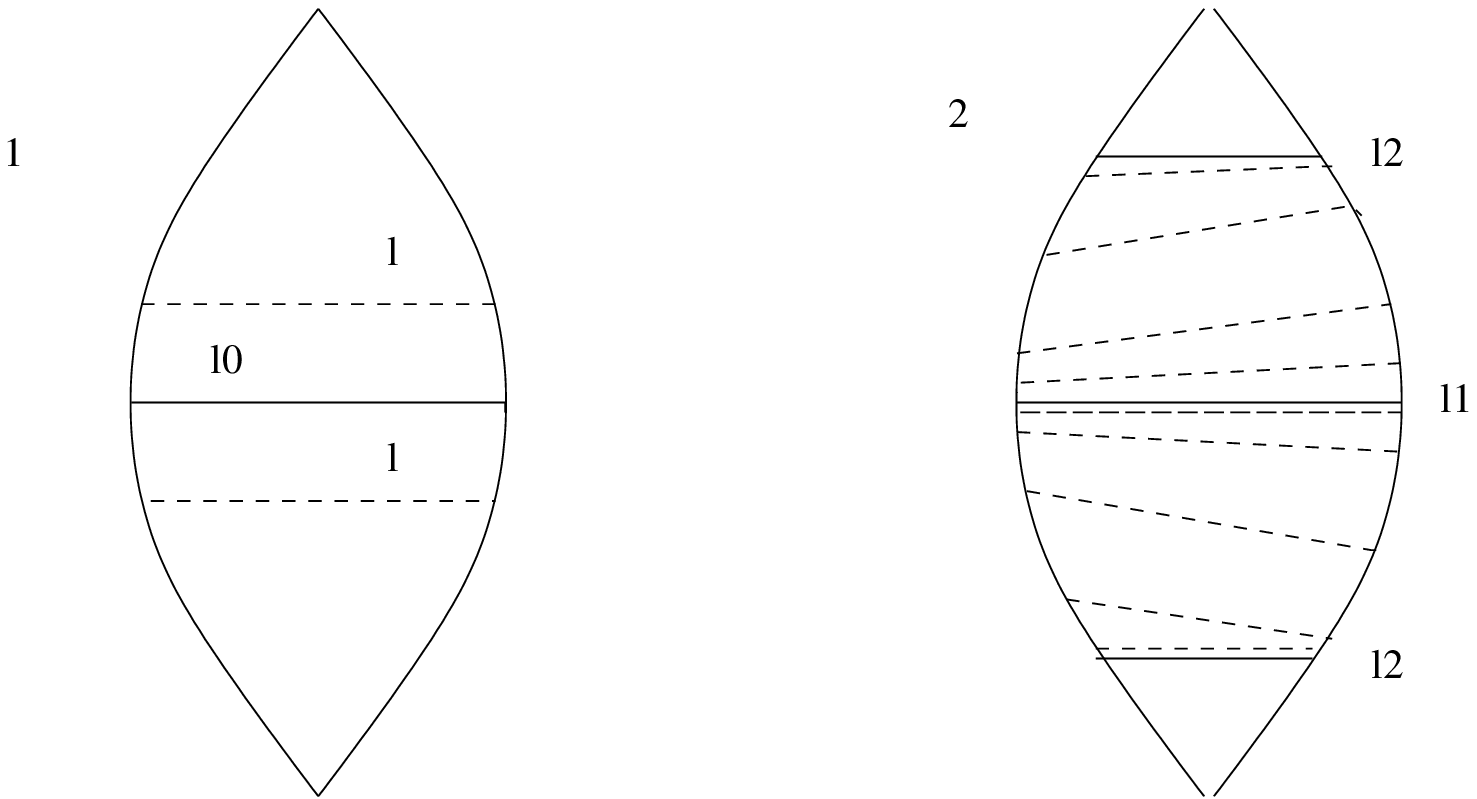}
        \caption{$\F_1,\F_2$ foliations on $\R P^2$: $Hol(L,\F_1)=\{e\}$, $Hol(L_0,\F_1)=\{e,g_0\}$, $g_0^2=e$, $Hol(L_1,\F_2)=\{e,g_1\}$, $g_1$ orientation reversing diffeomorphism, $Hol(L_2,\F_2)=\{e,g_2\}$, $g_2$ generator of unilateral holonomy.}
        \label{holonomy}
}
\end{center}
  }\end{minipage}%
  \begin{minipage}[t]{.1\linewidth}{\hspace{.1\linewidth}}\end{minipage}%
  \begin{minipage}[t]{.45\linewidth}
  {
\begin{center}
{        
        \psfrag{1}{$p_2$}
        \psfrag{2}{$p_3$}
        \psfrag{3}{$p_1$}
        \psfrag{4}{$q$}
        \psfrag{a}{$a$}
        \includegraphics[scale=.30]{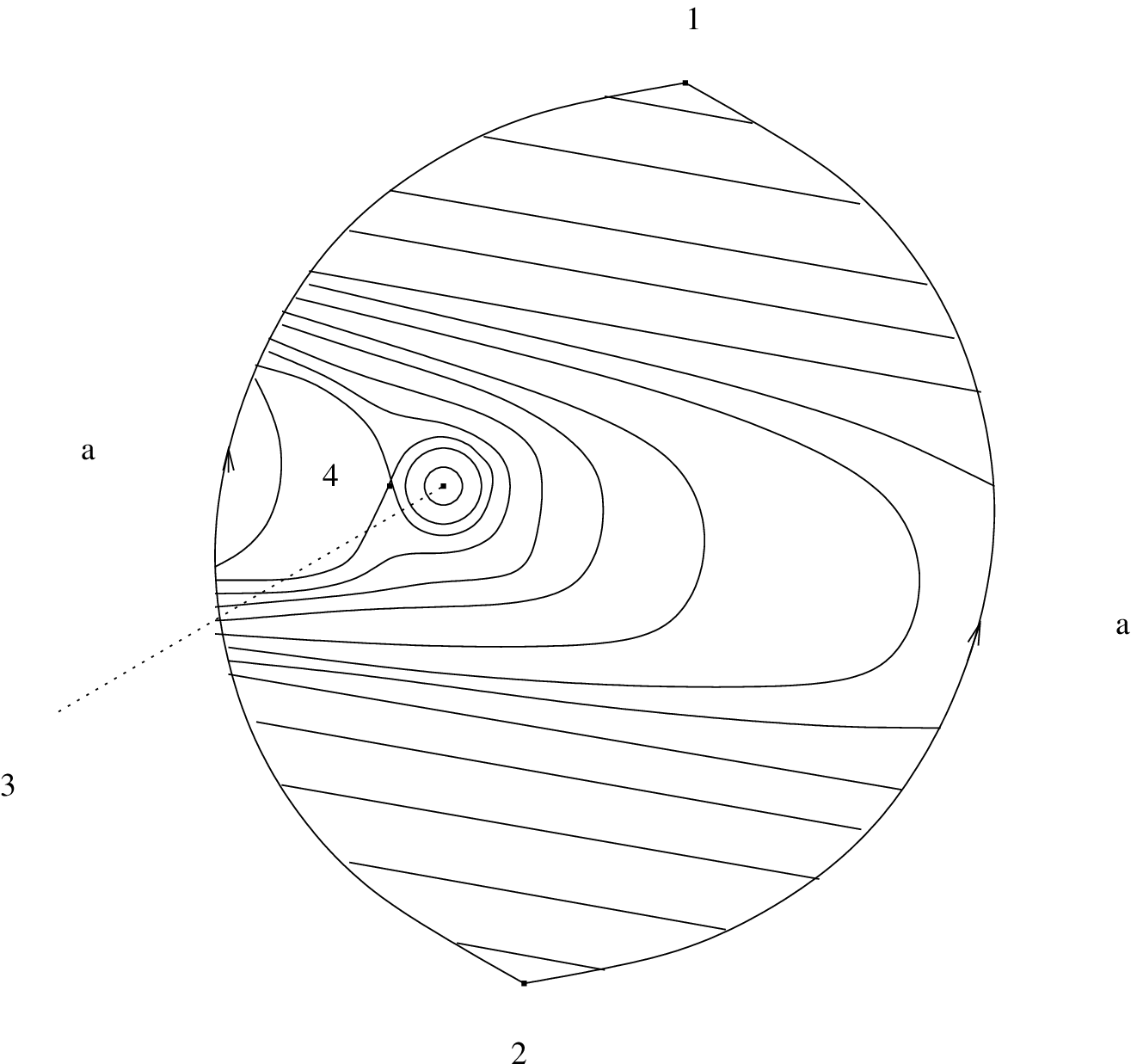}
        \caption{A singular foliation of the sphere $S^2$, which does not admit
         a first integral. With the same spirit, a singular foliation on $S^3$ may be given.}
 \label{less}
}
\end{center}  
}\end{minipage}
\end{figure}

It is well known that a basic tool in the study of foliations is the holonomy of a leaf (in the sense of Ehresmann). If $L$ is a leaf of a codimension $k$ foliation $(M,\F)$, the holonomy $Hol(L,\F)=\Phi(\pi_1(L))$, is the image of a representation, $\Phi:\pi_1(L) \rightarrow Germ(\R^k,0)$, of the fundamental group of $L$ into the germs of diffeomorphisms of $\R^k$, fixing the origin. Let $x \in L$ and $\Sigma _x$ be a section transverse to $L$ at $x$; with abuse of notation, we will write that a diffeomorphism $g: Dom(g) \subset \Sigma_x \rightarrow \Sigma_x$, fixing the origin, is an element of the holonomy group. For codimension one foliations ($k=1$), we may have: {\em{(i)}} $Hol(L,\F)=\{e \}$, {\em{(ii)}} $Hol(L,\F)=\{e,g \}$, with $g^2=e, g \neq e$, {\em{(iii)}} $Hol(L,\F)=\{e,g \}$, where $g^n \neq e$, $\forall n$, and $g$ is a (orientation preserving or reversing) diffeomorphism. In particular, among orientation preserving diffeomorphisms, we might find  a $g: \Sigma_x \rightarrow \Sigma_x$, such that $g$ is the identity on one component of $\Sigma_x \setminus \{x \}$ and it is not the identity on the other; in this case, we say that $L$ has {\em{unilateral holonomy}}  (see figure \ref{holonomy} for some examples). 
We recall Reeb Stability Theorems (cfr., for example, \cite{Cam} or \cite{Mo-Sca}).\\{\bf{Theorem \ref{uno}.1 (Reeb Local Stability)}} {\em{Let $\F$ be a $C^1$, codimension $k$ foliation of a manifold $M$ and $F$ a compact leaf with finite holonomy group. There exists a neighborhood $U$ of $F$, saturated in $\F$ (also called {\em{invariant}}), in which all the leaves are compact with finite holonomy groups. Further, we can define a retraction $\pi:U \rightarrow F$ such that, for every leaf $F' \subset U$, $\pi|_{F'} : F' \rightarrow F$ is a covering with a finite number of sheets and, for each $y \in F$, $\pi^{-1}(y)$ is homeomorphic to a disk of dimension $k$ and is transverse to $\F$. The neighborhood $U$ can be taken to be arbitrarily small.}}

The last statement means in particular that, in a neighborhood of the point corresponding to a compact leaf with finite holonomy, the space of leaves is Hausdorff.

 Under certain conditions the Reeb Local Stability Theorem may replace the Poincar\'e Bendixon Theorem \cite{Palis} in higher dimensions. This is the case of codimension one, singular foliations $(M^n,\F)$, with $n \geq 3$, and some center-type singularity in $Sing(\F)$.\\
{\bf{Theorem \ref{uno}.2 (Reeb Global Stability)}} {\em{Let $\F$ be a $C^1$, codimension one foliation of a closed manifold, $M$. If $\F$ contains a compact leaf $F$ with finite fundamental group, then all the leaves of $\F$ are compact, with finite fundamental group. If $\F$ is transversely orientable, then every leaf of $\F$ is diffeomorphic to $F$; $M$ is the total space of a fibration $f:M \rightarrow S^1$ over $S^1$, with fibre $F$, and $\F$ is the fibre foliation, $\{f^{-1}(\theta)| \theta \in S^1 \}$.}}

This theorem holds true even when $\F$ is a foliation of a manifold with boundary, which is, a priori, tangent on certain components of the boundary and transverse on other components \cite{God}. In this setting, let $\h ^l=\{(x^1, \dots , x^l) \in \R^l|x^l \geq 0 \}$. Taking into account definition \ref{inizio}.1, we say that a foliation of a manifold with boundary is {\em{tangent}}, respectively {\em{transverse}} {\em{to the boundary}}, if there exists a differentiable atlas $\{(U_i, \phi_i)\}_{i \in I}$, such that property (1) of the above mentioned definition holds for domains $U_i$ such that $U_i \cap \partial M = \varnothing$, while $\phi_i(U_i) = \bo ^{n-k} \times \h ^k$, respectively, $\phi_i(U_i) = \h ^{n-k} \times \bo ^k$ for domains such that $U_i \cap \partial M \neq \varnothing$. Moreover, we ask that the change of coordinates has still the form described in property (2). Recall that $\F|_{\partial M}$ is a regular codimension $k-1$ (respectively, $k$) foliation of the $(n-1)$-dimensional boundary. After this, it is immediate to write the definition for foliations which are tangent on certain components of the boundary and transverse on others.\\
Observe that, for foliations tangent to the boundary, we have to replace $S^1$ with $[0,1]$ in the second statement of the Reeb Theorem \ref{uno}.2 (see Lemma \ref{reeb}.6). 

We say that a component of $Sing(\F)$ is {\em{weakly stable}} if it admits a neighborhood, $U$, such that $\F|_{U}$ is a foliation with all leaves compact. The problem of global stability for a foliation with weakly stable singular components may be reduced to the case of foliations of manifolds with boundary, tangent to the boundary. It is enough to cut off an invariant neighborhood of each singular component.

 Holonomy is related to transverse orientability by the following:\\
{\bf{Proposition \ref{uno}.3}} {\em{Let $L$ be a leaf of a codimension one (Morse) foliation $(M,\F)$. If $Hol(L,\F)=\{e,g\}$, where $g^2=e$, $g \neq e$, then $\F$ is non-transversely orientable. Moreover, if $\pi:M \rightarrow
M/{\F}$ is the projection onto the space of leaves, then $\partial (M/{\F})
\neq \varnothing$ and $\pi(L) \in \partial (M/{\F})$}}.\\
{\em{Proof.}} We choose $x \in L$ and a segment $\Sigma_x$, transverse to the foliation at $x$. Then $g: \Sigma_x \rightarrow \Sigma_x$ turns out to be  
$g(y)= -y$. Let $y \rightarrow N_y$ a 1-plane field complementary to the tangent plane field $y \rightarrow T_y\F_y$. Suppose we may choose a vector field $y \rightarrow X(y)$ such that $N_y= \textrm{span} \{X(y) \}$. Then it shoud be $X(x)= -X(x)=(\textrm{d}g)_x(X(x))$, a contraddiction. Consider the space of leaves near $L$; this space is the quotient of $\Sigma_x$
with respect to the equivalence relation $\sim$ which identifies points on $\Sigma_x$ of the
same leaf. Then $\Sigma_x/_ \sim$ is a segment of type $(z,x]$ or $[x,z)$, where $\pi^{-1}(x)=L$.

\begin{figure}[t!]
\begin{minipage}[t]{.45\linewidth}
  {
\begin{center}
{
        \psfrag{b}{$b$}
        \psfrag{a}{$a$}
        \includegraphics[scale=.5]{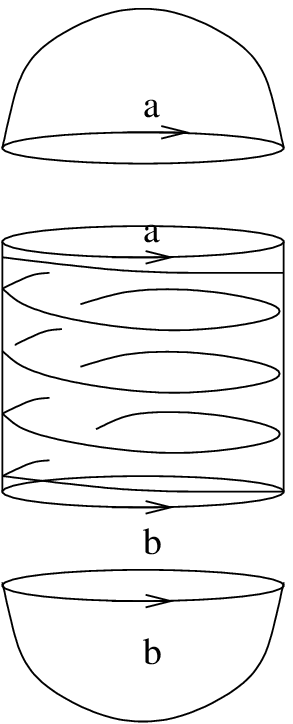}
        \caption{$n=2$: a singular foliation with center-type singularities,
        having no first integral.}
 \label{nofunction}
}
\end{center}
  }
  \end{minipage}%
  \begin{minipage}[t]{.1\linewidth}{\hspace{.1\linewidth}}\end{minipage}%
  \begin{minipage}[t]{.45\linewidth}
  {
\begin{center}
{ 
        \psfrag{1}{$R_1$}
        \psfrag{2}{$R_2$}
        \psfrag{3}{$R_3$}
        \psfrag{s}{$\W^s(q))$}
        \psfrag{i}{$\W^u(q)$} 
        \psfrag{p}{$p$}
        \psfrag{q}{$q$}
        \psfrag{u}{$U$} 
        \psfrag{l}{$L$} 
        \psfrag{f}{$F$}
        \includegraphics[scale=.35]{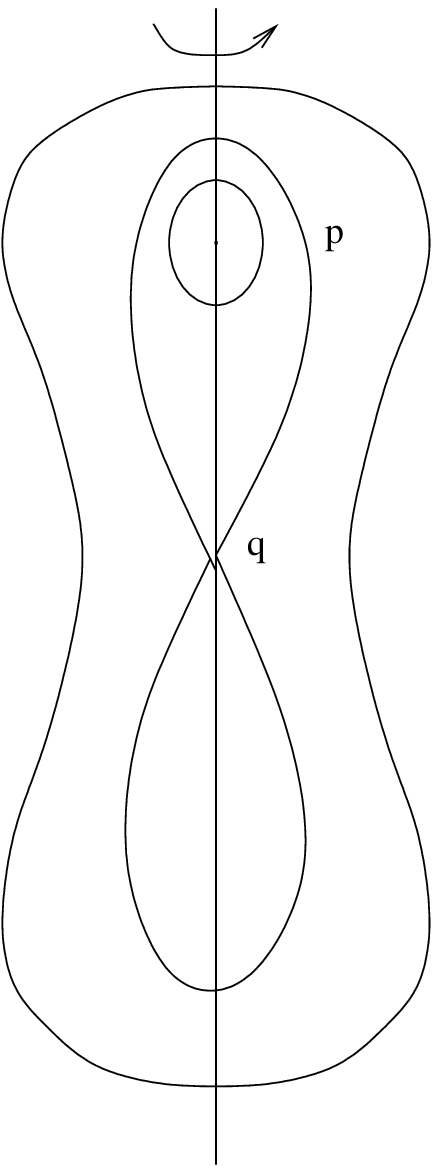}
        \caption{A trivial couple center-saddle $(p,q)$ (Theorem \ref{tre}.5, case  {\em{(i)}}).}
     \label{coupling1}
}
\end{center}
  }\end{minipage}
\end{figure}

At last we recall a classical result due to Reeb.\\
{\bf{Theorem \ref{uno}.4 (Reeb Sphere Theorem) \cite{Reeb}}} {\em{A transversely orientable Morse
foliation on a closed manifold, $M$, of dimension $n \geq 3$, having only
centers as singularities, is
homeomorphic to the $n$-sphere.}}\\
This result is proved by showing that the foliation considered must be given by a Morse function with
only two singular points, and therefore thesis follows by Morse theory.
Notice that the theorem still holds true for $n=2$, with a different proof. In particular, the foliation need not to be given by a function (see figure \ref{nofunction}).
\section{Arrangements of singularities}\label{tre}
In section \ref{quattro} we will study the elimination of singularities for Morse foliations. To this aim we will describe here how to identify special ``couples'' of singularities and we will study the topology of the neighbouring leaves.\\
{\bf{Definition \ref{tre}.1}} Let $n=\dim M, n \geq 2$. We define the set
$\C(\F)\subset M$ as the union of center-type singularities and  leaves
diffeomorphic to $S^{n-1}$ (with trivial holonomy if $n=2$) and for a center singularity, $p$, we
denote by  $\C_p(\F)$ the connected component of $\C(\F)$ that
contains $p$.\\
{\bf{Proposition \ref{tre}.2}} {\em{Let $\F$ be a Morse foliation on a manifold $M$. We have:\\
(1) $\C(\F)$ and $\C_p(\F)$ are open in $M$.\\
(2) $\C_p(\F) \cap \C_q(\F) \neq \varnothing$
if and only if $\C_p(\F)=\C_q(\F)$. $\C_p(\F)=M$ if and only if $\partial \C_p(\F)= \varnothing$. In this case the singularities of $\F$ are centers and the leaves are all diffeomorphic to $S^{n-1}.$\\
(3) If $q \in Sing(\F) \cap \partial \C_p(\F)$, then $q$ must be a saddle; in this case $\partial \C_p(\F) \cap Sing(\F)= \{ q \}$. Moreover, for $n \geq 3$ and $\F$ transversely orientable, $\partial \C_p(\F) \neq \varnothing$ if and only if $\partial \C_p(\F) \cap Sing(\F) \neq \varnothing$. In these hypotheses, $\partial \C_p(\F)$ contains at least one separatrix of the saddle $q$.\\
(4) $\partial \C_p(\F) \setminus \{q \}$ is closed in $M \setminus \{q \}$.}}\\
{\em{Proof.}} (1) $\C(\F)$ is open by the Reeb Local Stability Theorem \ref{uno}.1. (3) If non-empty, $\partial \C_p(\F) \cap Sing(\F)$ consists of a single saddle $q$, as there are no saddle connections. The second part follows by the Reeb Global Stability Theorem for manifolds with boundary and the third by the Morse Lemma. (4) By the Transverse Uniformity Theorem (see, for example, \cite{Cam}), it follows that the intrinsic topology of $\partial \C_p(\F) \setminus \{q \}$ coincides with its natural topology, as induced by $M \setminus \{q \}$.

We recall the following (cfr., for example \cite{Mo-Sca}):\\
{\bf{Lemma \ref{tre}.3 (Holonomy Lemma)}} {\em{Let $\F$ be a codimension one, transversely orientable foliation on $M$, let $A$ be a leaf of $\F$ and $K$ be a compact and path-connected set. If $g:K \rightarrow A$ is a $C^1$ map homotopic to a constant in $A$, then $g$ has a {\em{normal extension}} i.e. there exist $\epsilon >0$ and a $C^1$ map $G:K \times [0,\epsilon]
\rightarrow M$ such that $G_t(x)=G^x(t)=G(x,t)$ has the following properties: {\em{(i)}} $G_0(K)=g$, {\em{(ii)}} $G_t(K) \subset A(t)$ for some leaf $A(t)$ of $\F$ with $A(0)=A$, {\em{(iii)}} $\forall x \in K$ the curve $G^x([0, \epsilon])$ is normal to $\F$.}}\\
\begin{figure}
\begin{minipage}[t]{.45\linewidth}
  {
\begin{center}
{
 \psfrag{1}{$F_1$}
        \psfrag{2}{$F_2$}
        \psfrag{p}{$p$}
        \psfrag{q}{$q$}
        \psfrag{l}{$L$}
        \includegraphics[scale=.35]{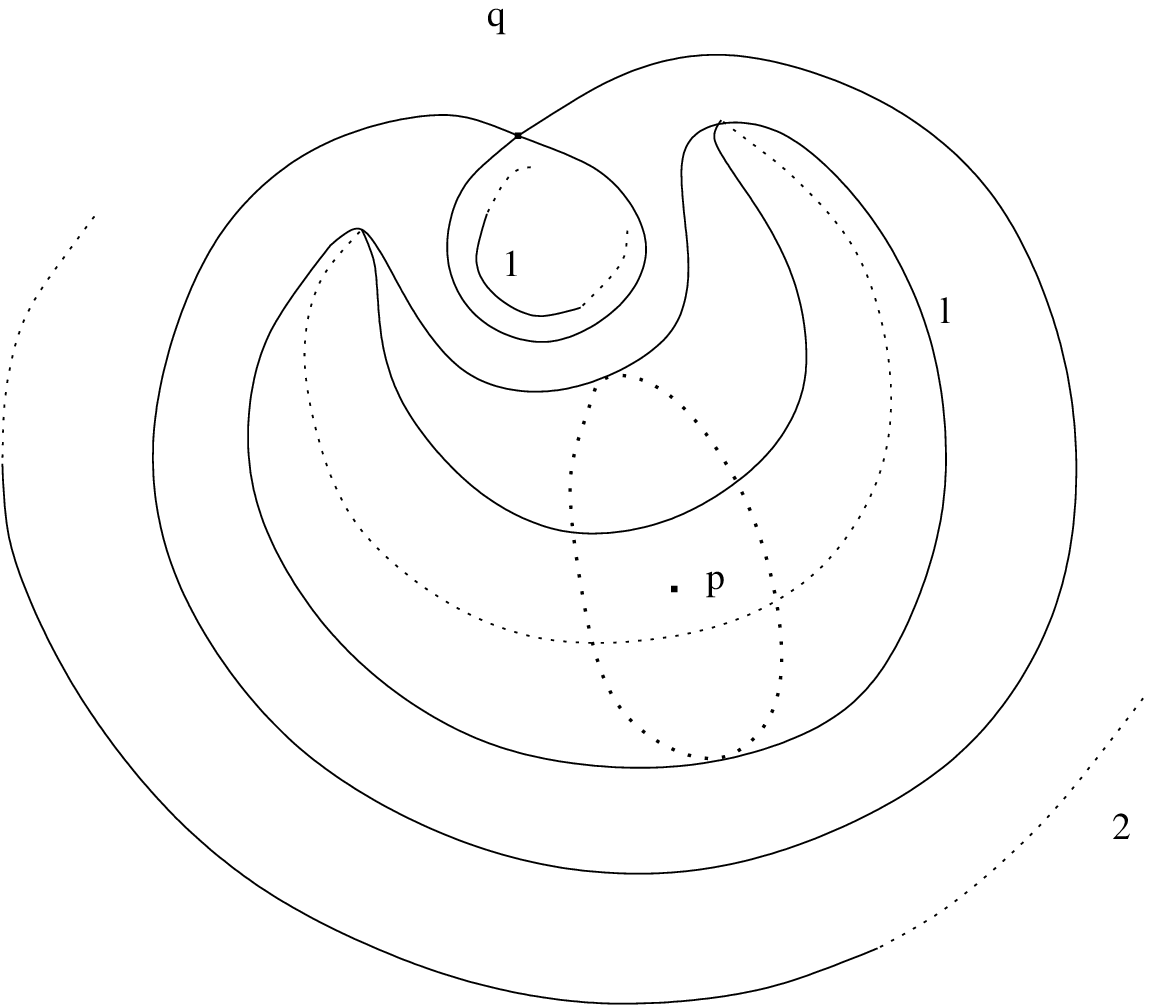}
        \caption{A saddle $q$ of index  $1$ ($n-1$), accumulating one center $p$ (Theorem \ref{tre}.5, case {\em{(ii)}}).}
        \label{coupling2}       
}\end{center}
  }\end{minipage}%
  \begin{minipage}[t]{.1\linewidth}{\hspace{.1\linewidth}}\end{minipage}%
  \begin{minipage}[t]{.45\linewidth}
  {
\begin{center}
{
        \psfrag{p}{$p$}
        \psfrag{q}{$q$}
        \psfrag{g}{$\gamma$}
        \includegraphics[scale=.25]{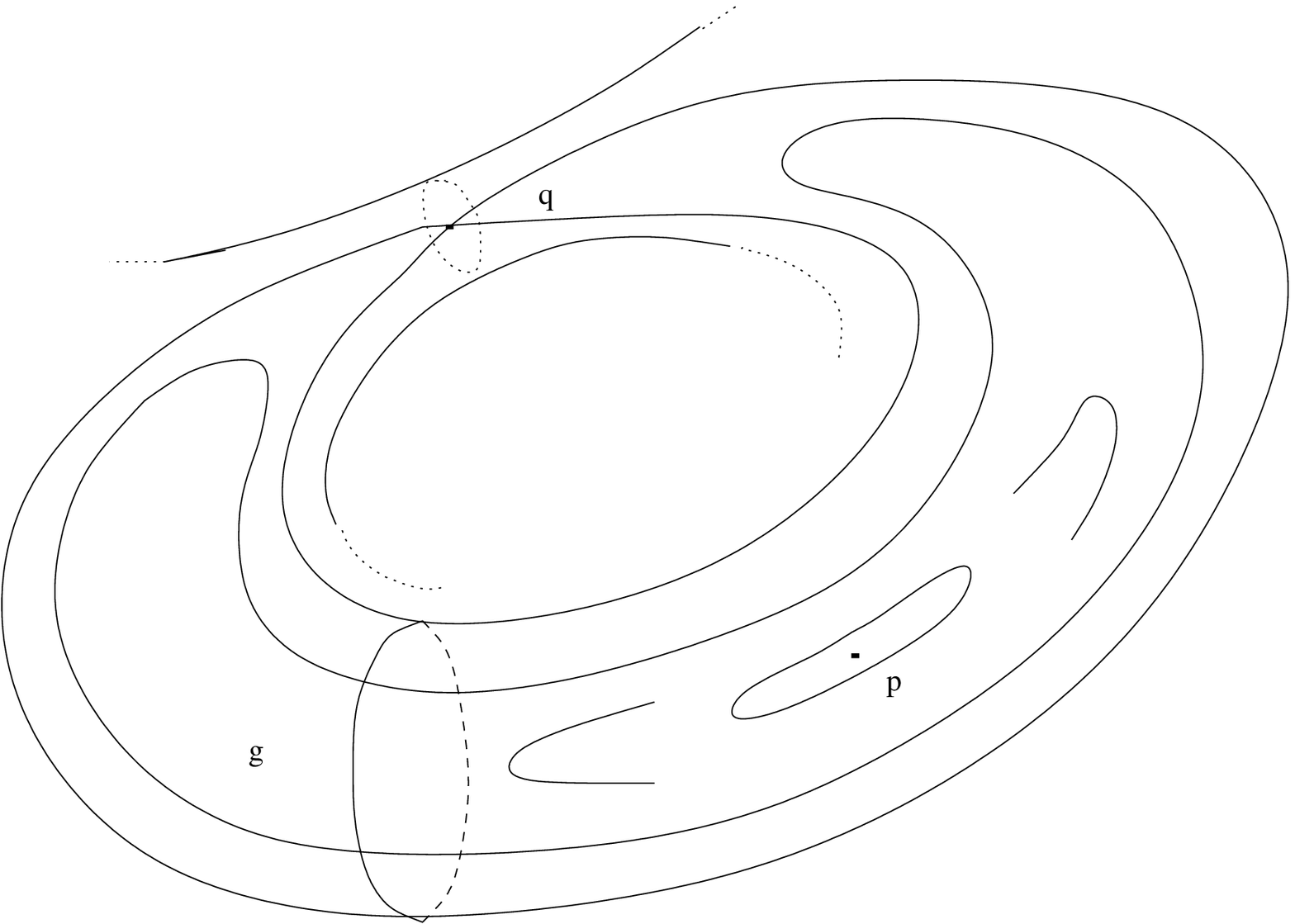}
        \caption{A saddle $q$ of index  $1$ ($n-1$), accumulating one center $p$ (Theorem \ref{tre}.5, case {\em{(iii)}}).}
        \label{coupling4}
}
\end{center}
  }\end{minipage}
\end{figure}
\begin{figure}[t!]
\begin{minipage}[t]{.45\linewidth}
  {
\begin{center}
{      
        \includegraphics[scale=.22]{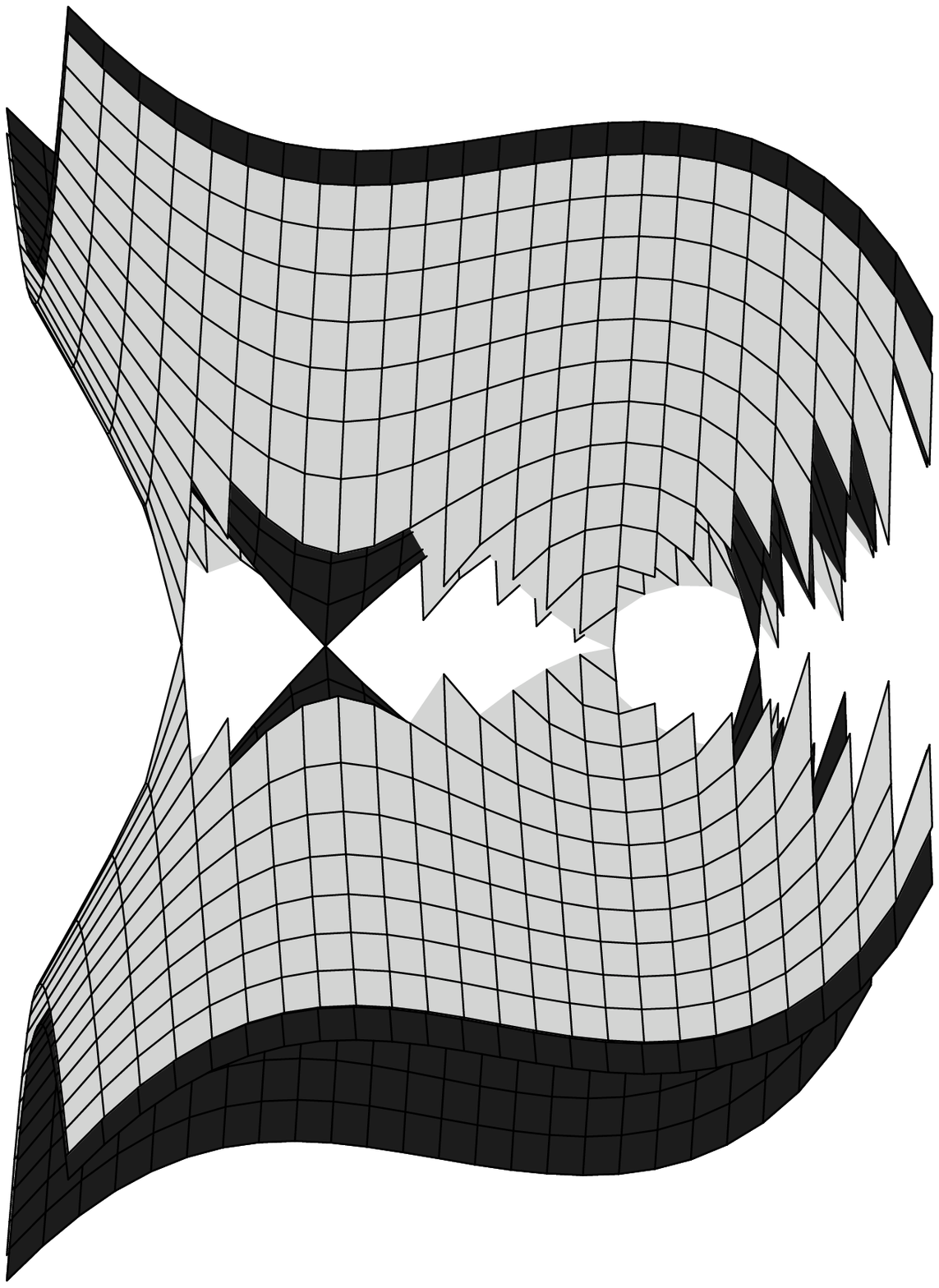}
        \caption{Two saddles in trivial coupling for the foliation defined
        by the function $f_\epsilon=- \frac{x^2}{2}+ \frac{y^3}{3}- \epsilon y+
        \frac{z^2}{2}$, ($\epsilon >0$).}
        \label{twosaddles}
}
\end{center}
  }
  \end{minipage}%
  \begin{minipage}[t]{.1\linewidth}{\hspace{.1\linewidth}}\end{minipage}%
  \begin{minipage}[t]{.45\linewidth}
  {
\begin{center}
{
        \psfrag{L_1}{$L_1$}
        \psfrag{L_2}{$L_2$}
        \psfrag{1}{$p$}
        \psfrag{2}{$q$}
        \psfrag{s_1}{$S_1$}
        \psfrag{s_2}{$S_2$}
        \psfrag{S}{$\Sigma$}
        \psfrag{n}{no intersection}
        \psfrag{L}{legenda}
        \includegraphics[scale=.55]{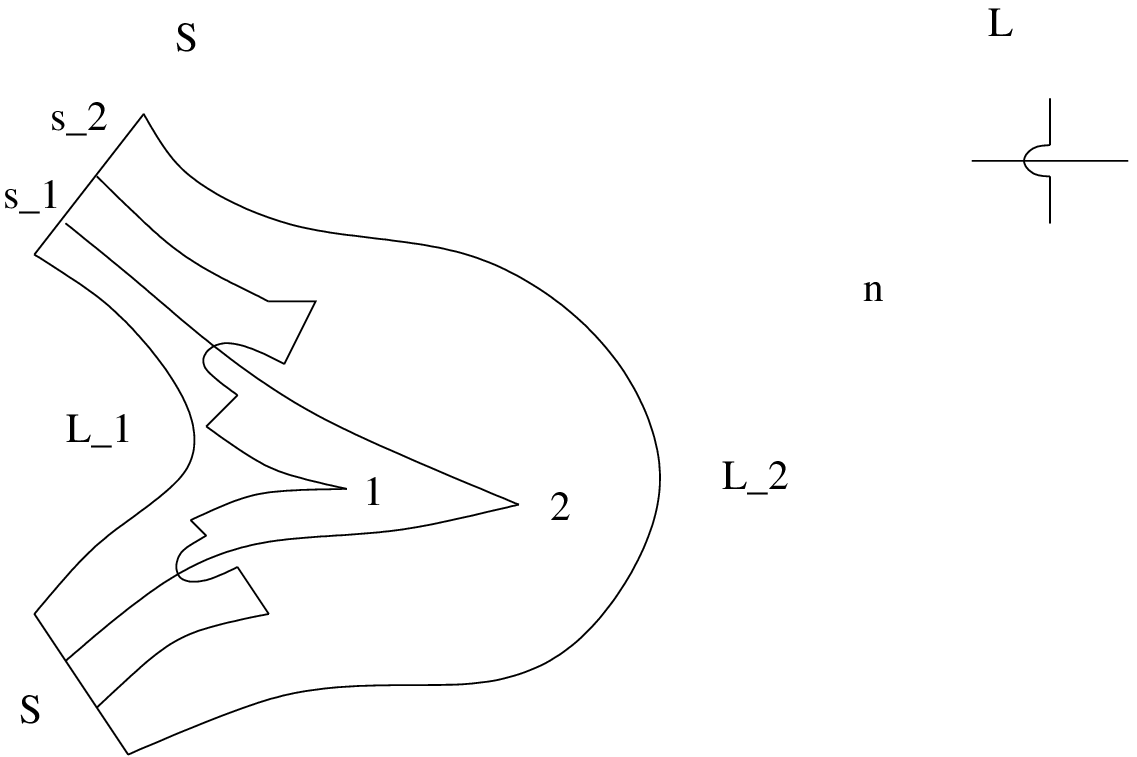}
        \caption{A dead branch of a trivial couple of saddles for a foliated manifold $(M^n,\F)$, $n \geq 3$.}
        \label{deadbranch}
}
\end{center}
  }\end{minipage}
\end{figure}
\hspace{3ex}For the case of center-saddle pairings we prove the following
descriptions of the separatrix:\\
{\bf{Theorem \ref{tre}.4}} {\em{Let $\F$ be a $C^ \infty$, codimension one, transversely
orientable, Morse foliation of a compact $n$-manifold, $M$, $n \geq
3$. Let $q$ be a saddle of index $l \notin
\{1, n-1 \}$, accumulating to one center
$p$. Let $L \subset \C_p(\F)$ be a spherical leaf intersecting a
neighborhood $U$ of $q$, defined by the Morse Lemma. Then $\partial \C_p(\F) \setminus \{q \}$ has a
single connected component (see figure \ref{LF}) and is homeomorphic to $S^{n-1}/S^{l-1}$. If $F$ is a leaf such that $F \cap \big (U \setminus \overline { \C_p(\F)} \big) \neq \varnothing$, then $F$ is homeomorphic to $\bo ^l \times S^{n-l-1} \cup_ \phi \bo ^l \times S^{n-l-1}$, where $\phi$ is a diffeomorphism of the boundary (for example, we may have $F \simeq S^l \times S^{n-l-1}$, but also $F \simeq S^{n-1}$, for $l=n/2$).}}\\
{\em{Proof.}} Let $\omega \in \Lambda^1(M)$ be a one-form defining the transversely orientable foliation. We choose a riemannian metric on $M$ and we consider the transverse vector field  $X_x=grad(\omega)_x$. We suppose $||X||=1$. In $U$, we have $X=h \cdot grad (f)$ for some real function $h>0$ defined on $U$. Further, we may suppose that $\partial U$ follows the orbits of $X$ in a neighborhood of $\partial \C_p(\F)$.\\
The Morse Lemma gives a local description of the foliation near its singularities; in particular the local topology of a leaf near a saddle of index $l$ is given by the connected components of the level sets of the function $f(x)=-x_1^2- \dots -x_l^2+x_{l+1}^2+ \dots +x_n^2$. If, for $c \geq 0$, we write $f^{-1}(c)=\{(x_1, \dots ,x_n) \in \R^n|x_1^2+ \dots +x_l^2+c=x_{l+1}^2+ \dots +x_n^2 \}$, it is easy to see that $f^{-1}(0)$ is homeomorphic to a cone over $S^{l-1} \times S^{n-l-1}$ and $f^{-1}(c) \simeq \bo ^l \times S^{n-l-1}$ ($c>0$). Similarly, we obtain $f^{-1}(c) \simeq \bo ^{n-l} \times S^{l-1}$ for $c<0$. Therefore, by our hypothesis on $l$, the level sets are connected; in particular the separatrix $S \supset f^{-1}(0)$ is unique and $\partial \C_p(\F)= S \cup \{ q \}$; moreover $U$ is splitted by $f^{-1}(0)$ in two different components. A priori, a leaf may intersect more than one component. As $\F$ is transversely orientable, the holonomy is an orientation preserving diffeomorphism, and then a leaf may intersect only non adiacent components; then this is not the case, in our hypotheses.\\
Let $L$ be a spherical leaf $\subset \C_p(\F)$ enough near $q$. Then $L \cap U \neq \varnothing$ and it is not restrictive to suppose it is given by $f^{-1}(c)$  for some $c<0$. We define the compact set $K=S^{n-1} \setminus \bo ^{n-l} \times S^{l-1} \simeq L \setminus U$. As $n \geq 3$, the composition $\xymatrix{
K \ar[rr]^ \simeq && L \setminus U \ar@{^{(}->}[rr]^\imath && L
}$ is homotopic to a constant in its leaf. By the proof of the Holonomy Lemma \ref{tre}.3, $L \setminus U$ projects diffeomorphically onto $A(\epsilon)=\partial \C_p(\F)$, by means of the constant-speed vector field, $X$. Together with the Morse Lemma, this gives a piecewise description of $\partial \C_p(\F)$, which is obtained by piecing pieces toghether. It comes out $\partial \C_p(\F) \simeq S^{n-1}/S^{l-1}$, a set with the homotopy type of $S^{n-1} \vee S^l$ (where $\vee$ is the wedge sum), simply connected in our hypotheses. Consequently, the map $K \times \{\epsilon \} \rightarrow  \partial \C_p(\F)$, obtained with the extension, admits, on turn, a normal extension. This completes the piecewise description of $F$. 

In case of presence of a saddle of index 1 or $n-1$, we have:\\
{\bf{Theorem \ref{tre}.5}} {\em{Let $\F$ be a $C^ \infty$, codimension one, transversely
orientable, Morse foliation of a compact $n$-manifold, $M$, $n
\geq 3$. Let $q$ be a saddle of index $1$ or $n-1$ accumulating to
one center $p$. Let $L \subset \C_p(\F)$ be a spherical leaf intersecting a
neighborhood $U$ of $q$, defined by the Morse Lemma. We may have: {\em{(i)}} $\partial \C_p(\F)$ contains a single separatrix of the
saddle (see figure \ref{coupling1}) and is homeomorphic to
$S^{n-1}$; {\em{(ii)}} $\partial \C_p(\F)$ contains both separatrices
$S_1$ and $S_2$ of the saddle (see figure \ref{coupling2}) and is
homeomorphic to $S^{n-1}/S^{n-2} \simeq S^{n-1} \vee S^{n-1}$. If this is the case, there exist two leaves
$F_i$ ($i=1,2$), such that $F_i$ and $L$ intersect different
components of $U \setminus S_i$ and we have that $F_i$ is
homeomorphic to $S^{n-1}$ ($i=1,2$); {\em{(iii)}} $q$ is a self-connected
saddle (see figure \ref{coupling4}) and $\partial \C_p(\F)$ is
homeomorphic to $S^{n-1}/S^0$. In this case we will refer to the couple $\Big(\overline{\C_p(\F)},\F|_{\overline{\C_p(\F)}}\Big)$ as a {\em{singular Reeb component}}. Moreover, $U \setminus \partial
\C_p(\F)$ has three connected components and $L$ intersects two of
them. If $F$ is a leaf intersecting the third
component of $U \setminus \partial \C_p(\F)$, then $F$ is
homeomorphic to $S^1 \times S^{n-2}$, or to $\R \times S^{n-2}$.}}\\
{\em{Proof.}} The proof is quite similar to the proof of the previous theorem. Nevertheless we give a brief sketch here. The three cases arise from the fact that $q$ has two local separatrices, $S_1$ and $S_2$, but not necessarily $\partial \C_p(\F)$ contains both of them. When this is the case, we may have that $S_1$ and $S_2$ belong to distinct leaves, or to the same leaf (in this case all spherical leaves contained in $\C_p(\F)$ intersect two different components of $U \setminus (S_1 \cup S_2)$ ). Using the Morse lemma, we construct the set $K$ for the application of the Holonomy Lemma \ref{tre}.3. We have, respectively: $K=\overline{\bo ^{n-1}}$, $K=K_1 \sqcup K_2= S^0 \times \overline{\bo ^{n-1}}$ (we apply twice the Holonomy Lemma), $K= \overline{\bo ^1} \times S^{n-2}$. In the first two cases, as $K$ is simply connected, the map $K \rightarrow L$, to be extended, is clearly homotopic to a constant in its leaf. Then $L \setminus U$ projects onto $\partial C_p(\F)$ and on neighbour leaves. This completes the piecewise description in case {\em{(i)}} and {\em{(ii)}}.\\In the third case, piecing pieces together after a first application of the Holonomy Lemma, we obtain $\partial \C_p(\F)\simeq S^{n-1}/S^0$ and $\partial \C_p(\F) \setminus \{ q \} \simeq \bo ^1 \times S^ {n-2}$, simply connected for $n \neq 3$. With a second application of the Holonomy Lemma ($n \neq 3$), $K$ projects diffeomorphically onto any neighbour leaf, $F$. The same also happens for $n=3$, because a curve $\gamma:S^1 \rightarrow \partial \C_p(\F)$, as the one depicted in figure \ref{coupling4}, is never a generator of the holonomy, which is locally trivial (a consequence of the Morse lemma). Nevertheless, there are essentially two ways to piece pieces together. We may have $F \simeq S^1 \times S^{n-2}$ or $F \simeq \R \times S^{n-2}$.

The last result gives the motivation for a new concept.\\
{\bf{Definition \ref{tre}.6}} In a codimension one singular foliation $\F$ it may happen that, for some leaf $L$ and  $q \in Sing(\F)$, the set $L \cup \{q \}$ is arcwise connected. Let $C=\{q \in Sing(\F)| L \cup \{q \}\textrm{ is arcwise connected} \}$. If for some leaf $L$ the set $C \neq \varnothing$, we define the corresponding {\em{singular leaf}} \cite{Wag} $S(L)= L \cup C$. In particular, if $\F$ is a transversely orientable Morse foliation, each singular leaf is given by $S(L)=L \cup \{ q \}$, for a single saddle-type singularity $q$, either selfconnected or not.

In the case of a transversely orientable Morse foliation $\F$ on $M$ ($n= \dim M \geq 3$), given a saddle $q$ and a separatrix $L$ of $q$, we may define a sort of holonomy map of the singular leaf $S(L)$. This is done in the following way.\\
As the foliation is Morse, in a neighborhood $U \subset M$ of $q$ there exists a (Morse) local first integral $f:U \rightarrow \R$, with $f(q)=0$. Keeping into account the structure of the level sets of the Morse function $f$ (see Theorem \ref{tre}.4 and Theorem \ref{tre}.5) we observe that there are at most three connected components in $U \setminus S(L)= U \setminus \{ f^{-1}(0)\}$ (notice that the number of components depends on the Morse index of $q$).\\
Let $\gamma: [0,1] \rightarrow S(L)$ be a $C^1$ path through the singularity $q$. At first, we consider  the case $\gamma([0,1]) \subset U$, $q= \gamma(t)$ for some $0<t<1$. For a point $x \in M \setminus Sing(\F)$, let $\Sigma_x$ be a transverse section at $x$. The set $\Sigma_x \setminus \{x \}$ is the union of two connected components, $\Sigma^+_x$ and $\Sigma ^-_x$ that we will denote by {\em{semi-transverse sections at $x$}}. For $x= \gamma(0) \in S(L)$ we have $f(x)=0$ and we can choose semi-transverse sections at $x$ in a way that $f(\Sigma^+_x)>0$ and $f(\Sigma^-_x)<0$. We repeat the construction for $y=\gamma(1)$, obtaining four semi-transverse sections, which are contained in (at most) three connected components of $U \setminus S(L)$. As a consequence, at least two of them are in the same component. By our choices, this happens for $\Sigma_x^-$ and $\Sigma_y^-$ (but we cannot exclude it happens also for $\Sigma_x^+$ and $\Sigma_y^+$). We define the {\em{semi-holonomy map}} $h^-:\Sigma^-_{\gamma(0)} \cup \gamma(0) \rightarrow \Sigma^-_{\gamma(1)} \cup \gamma(1)$ by setting $h^-(\gamma(0))=\gamma(1)$ and $h^-(z)=h(z)$ for $z \in \Sigma^-_{\gamma(0)}$, where $h:\Sigma^-_{\gamma(0)} \rightarrow \Sigma^-_{\gamma(1)}$ is a classic holonomy map (i.e. such that for a leaf $F$, it is $h(F \cap \Sigma^-_{\gamma(0)})=F \cap \Sigma^-_{\gamma(1)}$). In the same way, if it is the case, we define $h^+$.\\
Consider now any curve $\gamma: [0,1] \rightarrow S(L)$. As $\F$ is transversely orientable, the choice of a semi-transverse section for the curve $\gamma([0,1]) \cap U$, may be extended continuously on the rest of the curve, $\gamma([0,1]) \setminus U$; with this remark, we use classic holonomy outside $U$. To complete the definition, it is enough to say what a semi-transverse section at the saddle $q$ is. In this way we allow $q \in \gamma(\partial[0,1])$. To this aim, we use the orbits of the transverse vector field, $grad(f)$. By the property of gradient vector fields, there exist points $t,v$ such that $\alpha(t)=\omega(v)=q$. Let $\Sigma _q^+$ ($\Sigma _q^-$) be the negative (positive) semi-orbit through $t$ ($v$). Each of $\Sigma _q^+$ and $\Sigma _q^-$, transverse to the foliation and such that $\overline {\Sigma _q^+} \cap \overline {\Sigma _q^-}= \{ q \}$, is a {\em{semi-transverse section}} at the saddle $q$.

 In this way, the {\em{semi-holonomy of a singular leaf}} $Hol^+(S(L), \F)$ is a representation of the fundamental group $\pi_1(S(L))$ into the germs of diffeomorphisms of $\R _{\geq 0}$ fixing the origin, $Germ(\R _{\geq 0},0)$.

Now we consider the (most interesting) case of a selfconnected separatrix $S(L)=\partial \C_p(\F)$, with $\partial \C_p(\F)$ satisfying the description of Theorem \ref{tre}.5, case {\em{(iii)}}. The singular leaf $\partial \C_p(\F)$, homeomorphic to $S^{n-1}/S^0$, has the homotopy type of $S^{n-1} \vee S^1$. We have $Hol^+(\partial \C_p(\F),\F)=\{ e, h^-_\gamma \}$, where $\gamma$ is the non trivial generator of the homotopy, and $h^-_\gamma$ is a map with domain contained in the complement $\complement \C_p(\F)$. The two options $h^-_\gamma=e$, $h^-_\gamma \neq e$ give an explanation of the two possible results about the topology of the leaves near the selfconnected separatrix.
    
\section{Realization and elimination of pairings of singularities}\label{quattro}
Let us describe one of the key points in our work, i.e. the
elimination procedure, which allows us to delete pairs of
singularities in certain configurations, and, this way, to lead us
back to simple situations as in the Reeb Sphere Theorem (\ref{uno}.4). We need the following notion \cite{Ca-Sca}:\\
{\bf{Definition \ref{quattro}.1}} Let $\F$ be a codimension one foliation with isolated singularities on a manifold $M^n$. By a {\em{dead branch}} of $\F$ we mean a region $R \subset M$ diffeomorphic to the product $\bo ^{n-1} \times \bo^1$, whose boundary, $\partial R \approx \bo ^{n-1} \times S^0 \cup S^{n-2} \times \bo ^1$, is the union of two invariant components (pieces of leaves of $\F$, not necessarily distinct leaves in $\F$) and, respectively, of transverse sections, $\Sigma \approx \{t \} \times \bo^1$, $t \in S^{n-2}$.\\
Let $\Sigma_i, i=1,2$ be two transverse sections. Observe that the holonomy from $\Sigma_1 \rightarrow \Sigma_2$ is always trivial, in the sense of the Transverse Uniformity Theorem \cite{Cam}, even if $\Sigma_i \cap S(L) \neq \varnothing$ for some singular leaf $S(L)$. In this case we refer to the holonomy of the singular leaf, in the sense above.

\begin{figure}
\begin{minipage}[t]{.45\linewidth}
  {
\begin{center}
{
        \psfrag{p}{$p$}
        \psfrag{q}{$q$}
        \psfrag{r}{$r$}
        \psfrag{s}{$s$}
        \psfrag{v}{$V$}
        \includegraphics[scale=.38]{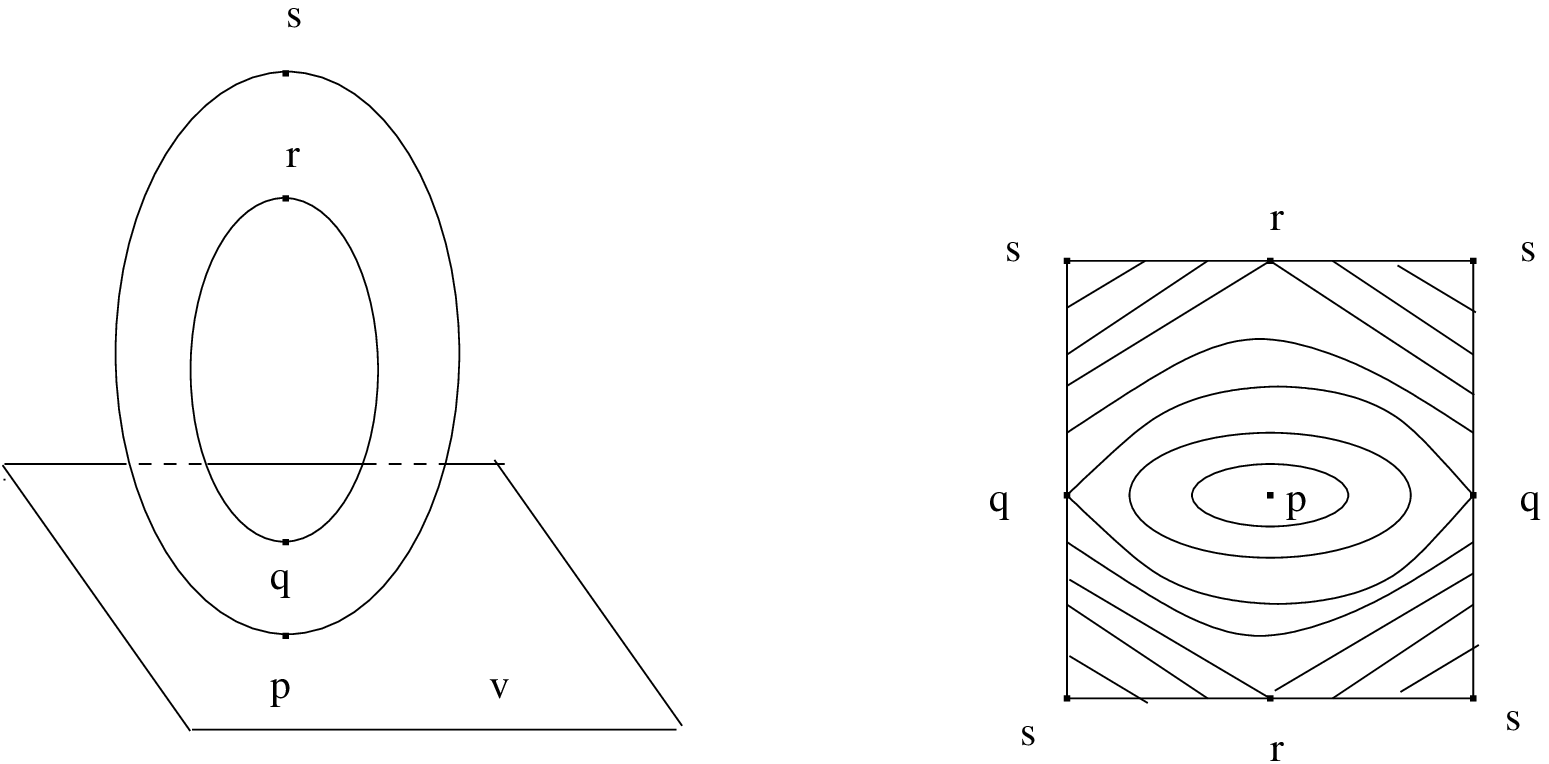}
        \caption{On the left: the height function on the plane V defines a foliation of the torus; on the right: a possible description of the foliation.}
        \label{milnor}
}      
\end{center}
  }\end{minipage}
  \begin{minipage}[t]{.1\linewidth}{\hspace{.1\linewidth}}\end{minipage}
  \begin{minipage}[t]{.45\linewidth}
  {      
\begin{center}
{
        \psfrag{1}{$P_1$}    
        \psfrag{2}{$P_2$}
        \psfrag{p}{$p$}
        \psfrag{q}{$q$}
        \psfrag{s}{$\partial \C_p(\F)$}
        \includegraphics[scale=.38]{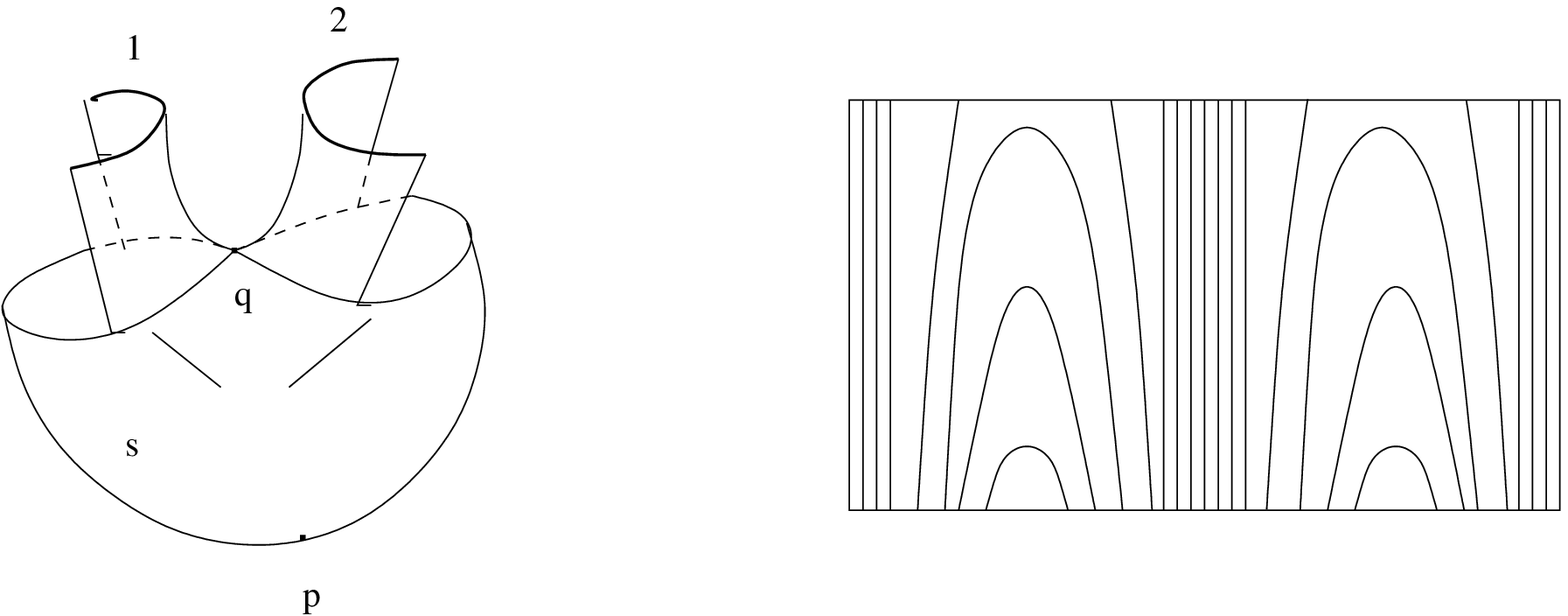}
        \caption{On the left, a dead branch for the selfconnected saddle $q$ of figure \ref{milnor}; on the right, the foliation obtained after the elimination of the two couples of singularities.}
        \label{milnor2}        
}    
\end{center}  
}\end{minipage}
\end{figure}

\begin{figure}
\begin{minipage}[t]{.45\linewidth}
  {
\begin{center}
{

 \includegraphics[scale=.55]{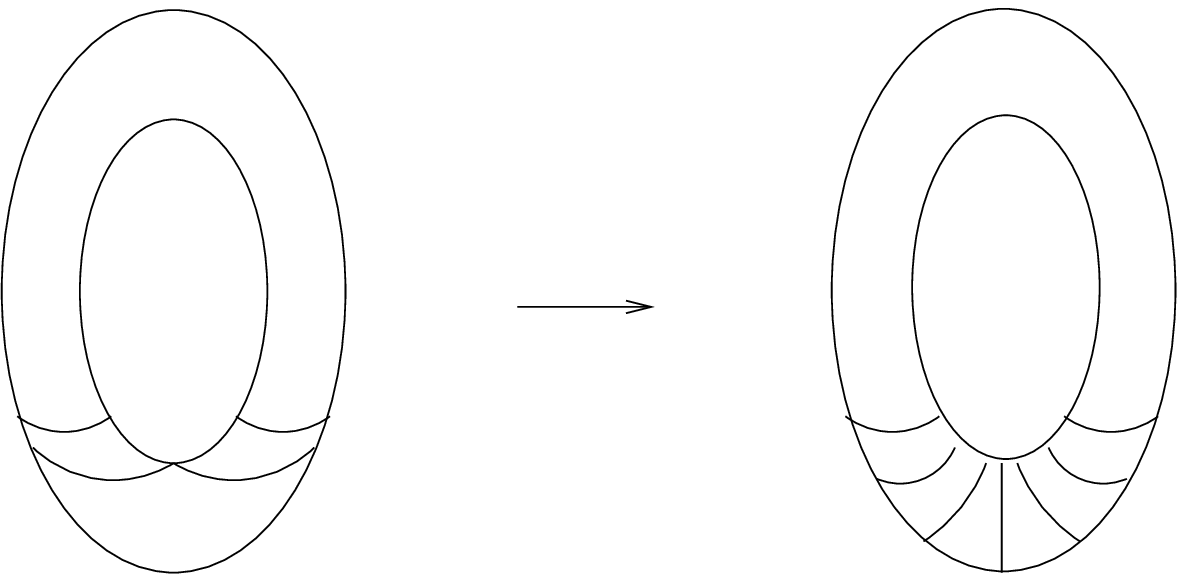}
        \caption{Elimination technique applied in case {\em{(ii)}} (Theorem \ref{tre}.5) for the foliation of figure \ref{milnor}.}
        \label{milnor4}      
}
\end{center}
  }\end{minipage}%
  \begin{minipage}[t]{.1\linewidth}{\hspace{.1\linewidth}}\end{minipage}%
  \begin{minipage}[t]{.45\linewidth}

  {
\begin{center}
{
        \psfrag{3}{$q$}
        \psfrag{1}{$p_1$}
        \psfrag{2}{$p_2$}
        \includegraphics[scale=.75]{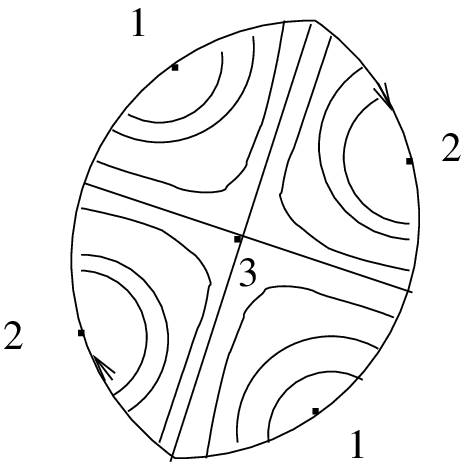}
        \caption{A foliation of $\R P 2$ with three singular points.}
        \label{rp2}
}
\end{center}  
}\end{minipage}
\end{figure}

 A first result includes known situations.\\
{\bf{Proposition \ref{quattro}.2}} {\em{Given a foliated manifold $(M^n,\F)$, with $\F$ Morse and transversely orientable, with $Sing(\F) \ni p,q$, where $p$ is a center and $q \in \partial \C_p(\F)$ is a saddle of index 1 or $n-1$, there exists a new foliated manifold $(M,\widetilde{\F})$, such that: {\em{(i)}} $\widetilde{\F}$ and $\F$ agree outside a suitable region $R$ of $M$, which contains the singularities $p,q$; {\em{(ii)}} $\widetilde{\F}$ is nonsingular in a neighborhood of $R$.}}\\
{\em{Proof.}} We are in the situations described by Theorem \ref{tre}.5. If we are in case {\em{(i)}}, the couple $(p,q)$ may be eliminated with the technique of the dead branch, as illustrated in \cite{Ca-Sca}. If we are in case {\em{(ii)}}, we observe that the two leaves $F_i, i=1,2$ bound a region, $A$, homeomorphic to an anulus, $S^{n-1} \times [0,1]$. We may now replace the singular foliation $\F|_A$ with the trivial foliation $\widetilde{\F}|_A$, given by $S^{n-1} \times \{t \}$, $t \in [0,1]$. If we are in case {\em{(iii)}}, we may replace the singular Reeb component with a regular one, in the spirit of \cite{Ca-Sca}. Even in this case, we may think the replacing takes place with the aid of a new sort of dead branch, the {\em{dead branch of the selfconnected saddle}}, that we describe with the picture of figure \ref{milnor2}, for the case of the foliation of the torus of figure \ref{milnor}, defined by the height Morse function \cite{Mil1}. Observe that the couples $(p,q)$ and $(r,s)$ of this foliation may be also seen as an example of the coupling described in Theorem \ref{tre}.5, case {\em{(ii)}}. In this case the elimination technique and the results are completely different (see figure \ref{milnor4}).\\    
{\bf{Definition \ref{quattro}.3}} If the couple $(p,q)$ satisfies the description of Theorem \ref{tre}.5, case {\em{(i)}} (and therefore may be eliminated with the technique of the dead branch), we will say that $(p,q)$ is a {\em{trivial couple}}.\\ 
\hspace{3ex}A new result is the construction of saddle-saddle
situations:\\
{\bf{Proposition \ref{quattro}.4}} {\em{Given a foliation $\F$ on an $n$-manifold $M^n$, there exists a
new foliation $\widetilde \F$ on $M$, with $Sing(\widetilde
\F)=Sing(\F) \cup \{p,q\}$, where $p$ and $q$ are a couple of
saddles of consecutive indices, {\em{connecting transversely}}
(i.e. such that the stable manifold of $p$, $\W^s(p)$, intersects
transversely the unstable manifold of $q$, $\W^u(q)$).}}\\
{\em{Proof.}} We choose the domain of (any) foliated chart, $(U,\phi)$. Observe that $R'=U$ ($\simeq \phi(U)$) is a dead branch for a foliation $\F_{\epsilon '}$, given (up to diffeomorphisms) by the submersion  $f_\epsilon= - x_1^2/2- \dots - x_{k-1}^2/2+(x_k^3/3- \epsilon x_k)+ x_{k+1}/2+ \dots + x_n^2/2$, for some $\epsilon =\epsilon '<0$.  We consider $\F_{\epsilon ''}$, given by taking $\epsilon=\epsilon '' >0$ in $f_\epsilon$, which presents a couple of saddles of consecutive indices, and we choose a dead branch $R''$ around them. We also choose a homeomorphism between $R'$ and $R''$ which sends invariant sets of $\F_{\epsilon '}$ into invariant sets of $\F_{\epsilon ''}$ in a neighborhood of the boundary. With a surgery, we may replace $\F_{\epsilon '}$ with $\F_{\epsilon ''}$.

The converse of the above poposition is preceded by the following\\
{\bf{Remark \ref{quattro}.5}} Given a foliation $\F$ on $M^n$ with two complementary saddle
singularities $p,q \in Sing(\F)$, having a strong stable connection
$\gamma$, there exist a neighborhood $U$ of $ p, q$ and $\gamma$
in $M^n$, a $\delta \in \R^+$ and a coordinate system $\phi: U \rightarrow \R^n$ taking
$p$ onto $(0,\dots, \phi^k=-\delta ,\dots,0)$, $q$ onto $(0,\dots, \phi^k=\delta ,\dots,0)$,
$\gamma$ onto the $x_k$-axis, $\{x_l=0\}_{l \neq k}$, and such that: {\em{(i)}} 
the stable manifold of $p$ is tangent to
$\phi^{-1}(\{x_l=0\}_{l>k})$ at $p$, {\em{(ii)}} the unstable manifold of
$q$ is tangent to $\phi^{-1}(\{x_l=0\}_{l<k})$ at $q$ (we are led to the situation considered in \cite{Mil2}, A first cancelation theorem). So
using the chart $\phi:U \rightarrow \R^n$ we may assume that we
are on a dead branch of $\R^n$ and the foliation $\F|_U$ is defined by $f_ \epsilon$,
for $\epsilon =\delta^2$. In this way the vector field $grad(f_
\epsilon)$ defines a transverse orientation in $U$. For a suitable $\mu>0$, the points $r_1=(0,\dots,\phi ^k=-\delta-\mu,
\dots,0)$ and $r_2=(0,\dots,\phi ^k=\delta+\mu, \dots,0)$ are such that the modification takes place in a region of $U$ delimited by $L_{r_i}, i=1,2$.\\
{\bf{Proposition \ref{quattro}.6}} {\em{Given a foliation $\F$ on $M^n$ with a couple of saddles  $p,q$ of complementary indices, having a strong stable connection, there exists a dead branch of the couple of saddles, $R \subset M$ and we can obtain a foliation $\widetilde{\F}$ on $M$ such that: {\em{(i)}} $\widetilde{\F}$ and $\F$ agree on $M \setminus R$; {\em{(ii)}} $\widetilde{\F}$ is nonsingular in a neighborhood of $R$; indeed $\widetilde{\F}|_R$ is conjugated to a trivial fibration; {\em{(iii)}} the holonomy of $\widetilde{\F}$ is conjugate to the holonomy of $\F$ in the following sense: given any leaf $L$ of $\F$ such that $L \cap (M \setminus R) \neq \varnothing$, then the corresponding leaf $\widetilde{L}$ of $\widetilde{\F}$ is such that $Hol( \widetilde{L},\widetilde{\F})$ is conjugate to $Hol(L,\F)$.}}\\
{\bf{Example \ref{quattro}.7 (Trivial Coupling of Saddles)}} Let $M=S^n, n \geq 3$. For $l=1, \dots ,n-2$ we may find a Morse foliation of $M=S^n$, invariant for the splitting $S^n=\overline{\bo ^{n-l}} \times S^l \cup_\phi S^{n-l-1} \times \overline{\bo ^{l+1}}$, where $\phi$ is a diffeomorphism of the boundary. In fact, by theorem \ref{tre}.4 or \ref{tre}.5, case {\em{(iii)}}, $\overline{\bo ^{n-l}} \times S^l$ admits a foliation with one center and one saddle of index $l$. Similarly, $S^{n-l-1} \times \bo ^{l+1}$ admits a foliation with a saddle of index $n-l-1$, actually a saddle of index $l+1$, after the attachment. We may eliminate the trivial couple of saddles and we are led to the well-known foliation of $S^n$, with a couple of centers and spherical leaves.\\ 
{\bf{Remark \ref{quattro}.8}} The elimination of saddles of consecutive indices is actually a generalization of the elimination of couples center-saddle, $(p,q)$ with $q \in \partial \C_p(\F)$. Indeed, we may eliminate $(p,q)$ only when the saddle $q$ has index $1$ or $n-1$. This means the singularities of the couple must have consecutive indices and, as $q \in \partial \C_p(\F)$, there exists an orbit of the transverse vector field having $p$ as $\alpha$-limit (backward) and $q$ as $\omega$-limit (forward), or viceversa. Such an orbit is a strong stable connection.
\section{Reeb-type theorems}\label{reeb}
We shall now describe how to apply our techniques to obtain some generalizations
of the Reeb Sphere Theorem (\ref{uno}.4) for the case of Morse foliations admitting both centers and saddles.\\
A first generalization is based on the following notion:\\
{\bf{Definition \ref{reeb}.1}} We say that an isolated singularity, $p$, of a $C^ \infty$,
codimension one foliation $\F$ on $M$ is a {\em{stable singularity}},
if there exists a neighborhood $U$ of $p$ in $M$ and a $C^ \infty$ function,
$f:U \rightarrow \R$, defining the foliation in $U$, such that $f(p)=0$ and
$f^{-1}(a)$ is compact, for $|a|$ small.
The following characterization of stable singularities can be found in \cite{Ca-Sca}.\\
{\bf{Lemma \ref{reeb}.2}} {\em{An isolated singularity $p$ of a function $f:U \subset \R^n
\rightarrow \R$ defines a stable singularity for $\textrm{d}f$, if and only
if there exists a neighborhood $V \subset U$ of $p$, such that,
$\forall x \in V$, we have either $\omega (x) = \{p \}$ or $\alpha
(x) = \{p \}$, where $\omega (x)$ (respectively  $\alpha (x)$) is the
$\omega$-limit (respectively $\alpha$-limit) of the orbit of the vector
field $grad (f)$ through the point $x$.}}

In particular it follows the well-known:\\
{\bf{Lemma \ref{reeb}.3}} {\em{If a function $f:U \subset \R^n \rightarrow \R$ has an isolated local maximum or minimum at $p \in U$ then $p$ is a stable singularity for $df$.}}

The converse is also true:\\
{\bf{Lemma \ref{reeb}.4}} {\em{If $p$ is a stable singularity, defined by the function $f$, then $p$ is a point of local maximum or minimum for $f$.}}\\
{\em{Proof.}} It follows immediately by Lemma \ref{reeb}.2 and by the fact that $f$ is monotonous, strictly increasing, along the orbits of $grad (f)$.

With this notion, we obtain\\
{\bf{Lemma \ref{reeb}.5}}  {\em{Let $\F$ be a codimension one, singular foliation on a manifold $M^n$. In a neighborhood of a stable singularity, the leaves of $\F$ are diffeomorphic to spheres.}}\\
{\em{Proof.}} Let $p \in Sing(\F)$ be a stable singularity.  By Lemma \ref{reeb}.4, we may suppose $p$ is a minimum (otherwise we use $-f$). Using a local chart around $p$, we may suppose we are on $\R^n$ and we may write the Taylor-Lagrange expansion around $p$ for an approximation of the function $f:U \rightarrow \R$ at the second order. We have $f(p+h)=f(p)+1/2 \langle h,H(p+ \theta h)h \rangle,$
where $H$ is the Hessian of $f$ and $0< \theta <1$. It follows $\langle h,H(p+ \theta h)h \rangle \geq 0$ in $U$. Then $f$ is convex and hence the sublevels, $f^{-1}(c)$, are also convex.\\
We consider the flow $\phi: \mathscr{D}(\phi) \subset \R \times U \rightarrow U$ of the vector field $grad(f)$. By the properties of gradient vector fields, in our hypothesis,  $\mathscr{D}(\phi) \supset (- \infty,0] \times U$ and $\forall x \in U$ there exists the $\alpha$-limit, $\alpha(x)=p$. For any $x \in f^{-1}(c)$, the tangent space, $T_x f^{-1}(c)$, to the sublevels of $f$ does not contain the radial direction, $\overrightarrow{px}$. This is obvious otherwise, for the convexity of $f^{-1}(c)$, the singularity $p$ should lie on the sublevel $f^{-1}(c)$, a contraddiction because, in this case, $p$ should be a saddle. Equivalently, the orbits of the vector field $grad(f)$ are transverse to spheres centered at $p$. An application of the implicit function theorem shows the existence of a smooth function $x \rightarrow t_x$, that assigns to each point $x \in f^{-1}(c)$ the (negative) time at which $\phi(t,x)$ intersects $S^{n-1}(p, \epsilon)$, where $\epsilon$ is small enough to have $\textrm{B}^{n}(p, \epsilon) \subsetneq R(f^{-1}(c))$, the compact region bounded by $f^{-1}(c)$ . The diffeomorphism between the leaf $f^{-1}(c)$ and the sphere $S^{n-1}(p, \epsilon)$ is given by the composition $x \rightarrow \phi(t_x,x)$. The lemma is proved.\\
{\bf{Lemma \ref{reeb}.6}} {\em{Let $\F$ be a codimension one, transversely orientable foliation of $M$, with all leaves closed, $\pi:M \rightarrow M/\F$ the projection onto the space of leaves. Then we may choose a foliated atlas on $M$ and a differentiable structure on $M/ \F$, such that $M/ \F$ is a codimension one compact manifold, locally diffeomorphic to the space of plaques, and $\pi$ is a $C^\infty$ map.}}\\
{\em{Proof.}} At first we notice that the space of leaves $M/\F$ (with the quotient topology) is a one-dimensional Hausorff topological space, as a consequence of the Reeb Local Stability Theorem \ref{uno}.1. As all leaves are closed and with no holonomy, we may choose a foliated atlas $\{(U_i, \phi_i)\}$ such that, for each leaf $L \in \F, L \cap U_i$ consists, at most, of a single plaque. Let $\pi:M \rightarrow M/ \F$ be the projection onto the space of leaves and $\pi_i:U_i \rightarrow \R$ the projection onto the space of plaques. With abuse of notation, we may write $\pi_i=p_2 \circ \phi_i$, where $p_2$ is the projection on the second component. As there is a 1-1 correspondence between the quotient spaces $\pi|_{U_i}(U_i)$ and $\pi_i(U_i)$, then, are homeomorphic. Let $V \subset M/ \F$ be open. The set $\pi^{-1}(V)$ is an invariant open set. We may find a local chart $(U_i,\phi_i)$ such that $\pi(U_i)=V$. We say that $(V, \pi_i \circ (\pi|_{U_i})^{-1})$ is a chart for the differentiable atlas with the required property. To see this, it is enough to prove that, if $(V, \pi_j \circ (\pi|_{U_j})^{-1})$ is another chart with the same domain, $V$, there exists a diffeomorphism between the two images of $V$, i.e. between $\pi_i \circ (\pi|_{U_i})^{-1}(V)$ and $\pi_j \circ (\pi|_{U_j})^{-1}(V)$. This is not obvious when $U_i \cap U_j= \varnothing$. Indeed, the searched diffeomorphism exists, and it is given by the Transverse Uniformity Theorem \cite{Cam}. Observe that, in coordinates, $\pi$ coincides with the projection on the second factor.\\
{\bf{Lemma \ref{reeb}.7}} {\em{Let $n \geq 2$. A weakly stable singularity for a foliation $(M^n, \F)$ is a stable singularity.}}\\
{\em{Proof.}} Let $p$ be a weakly stable singularity, $U$ a neighborhood of $p$ with all leaves compact. We need a local first integral near $p$. As a consequence of the Reeb Local Stability Theorem \ref{uno}.1, we can find an (invariant) open neighborhood $V \subset U$ of $p$, whose leaves have all trivial holonomy. The set $V \setminus \{ p \}$ is open in $M^*=M \setminus Sing(\F)$. Let $\F^*= \F \setminus Sing(\F)$; the projection $\pi^*:M^* \rightarrow M^*/\F^*$ is an open map (see, for example \cite{Cam}). As a consequence of Lemma \ref{reeb}.6, the connected (as $n \geq 2$) and open set $\pi^* (V \setminus \{ p \})$ is a $1$-dimensional manifold with boundary, i.e. it turns out to be an interval, for example $(0,1)$. Now, we extend smoothly  $\pi^*$ to a map $\pi$ on $U$. In particular, let $W \subsetneq V$ be a neighborhood of $p$. If (for example) $\pi^*(W \setminus \{p \})=(0,b)$ for some $b<1$, we set $\pi(p)=0$. Thesis follows by lemma \ref{reeb}.3.\\
{\bf{Theorem \ref{reeb}.8}} {\em{Let $M^n$ be a closed $n$-dimensional manifold, $n \geq
3$. Suppose that $M$ supports a $C^ \infty$, codimension one,
transversely orientable foliation, $\F$, with non-empty singular
set, whose elements are, all, weakly stable singularities. Then $M$ is
homeomorphic to the sphere, $S^n$.}}\\
{\em{Proof.}} By hypothesis, $\forall p \in Sing(\F)$, $p$ is a weakly stable singularity. Then it is a stable singularity. By lemma \ref{reeb}.5, in an invariant neighborhood $U_p$ of $p$, the leaves are diffeomorphic to spheres. Now we can proceed as in the proof of the Reeb Sphere Theorem \ref{uno}.4.\\ 
{\bf{Theorem \ref{reeb}.9 (Classification of codimension one foliations with all leaves compact)}} {\em{Let $\F$ be a (possibly singular, with isolated singularities) codimension one foliation of $M$, with all leaves compact. Then all possible singularities are stable. If $\F$ is (non) transversely orientable, the space of leaves is (homeomorphic to $[0,1]$) diffeomorphic to $[0,1]$ or $S^1$. In particular, this latter case ocurs if and only if $\partial M, Sing(\F)= \varnothing$. In all the other cases, denoting by $\pi:M \rightarrow [0,1]$ the projection onto the space of leaves, it is $Hol(\pi^{-1}(x), \F)=\{e \}, \forall x \in (0,1)$. Moreover, if $x=0,1$, we may have: {\em{(i)}} $\pi^{-1}(x) \subset \partial M \neq \varnothing$ and $Hol(\pi^{-1}(x), \F)=\{e \}$; {\em{(ii)}} $\pi^{-1}(x)$ is a (stable) singularity; {\em{(iii)}} $Hol(\pi^{-1}(x),\F)=\{e,g \}$, $g\neq e, g^2=e$ (in this case, $\forall y \in (0,1)$, the leaf $\pi^{-1}(y)$ is a two-sheeted covering of $\pi^{-1}(x)$.}}\\
{\em{Proof.}} If $\F$ is transversely orientable, by the Reeb Global Stability Theorem \ref{uno}.2 and Lemma \ref{reeb}.6, the space of leaves is either diffeomorphic to $S^1$ or to $[0,1]$. In particular, $M/ \F \approx S^1$ if and only if $M$ is closed and $\F$ non singular. When this is not the case, $M/ \F \approx [0,1]$, and there are exactly two points ($\partial [0,1]$) which come from a singular point and/or from a leaf of the boundary.\\
If $\F$ is non transversely orientable, there is at least one leaf with (finite) non trivial holonomy, which corresponds a boundary point in $M/ \F$ to (by Proposition \ref{uno}.3). By the proof of Lemma \ref{reeb}.6, the projection is not differentiable and the space of leaves $M/ \F$, a Hausdorff topological $1$-dimensional space, turns out to be an orbifold (see \cite{Thu}). We pass to the transversely orientable double covering, $p: (\widetilde{M}, \widetilde{\F}) \rightarrow (M,\F)$. The foliation $\widetilde{\F}$, pull-back of $\F$, has all leaves compact, and singular set empty or with stable components; therefore we apply the first part of the classification to $\widetilde{M}/\widetilde{\F}$. Both if $\widetilde{M}/\widetilde{\F}$ is diffeomorphic to $S^1$ or to $[0,1]$, $M/ \F$ is homeomorphic to $[0,1]$, but (clearly) with different orbifold structures.      
\hspace{3ex}

Before going on with our main generalization of the Reeb Sphere Theorem \ref{uno}.4, which extends a similar result of Camacho and Sc\'ardua \cite{Ca-Sca} concerning the case $n=3$, we need to recall another result, that we are going to generalize.\\
As we know, the Reeb Sphere Theorem, in its original statement, consideres the effects (on the topology of a manifold $M$) determined by the existence, on $M$, of a real valued function with exactly two non-degenerate singular points. A very similar problem was studied by Eells and Kuiper \cite{Ku-Ee}. They considered manifolds admitting a real valued function with exactly three non-degenerate singular points.They obtained very interesting results. Among other things, it sticks out the obstruction they found about the dimension of $M$, which must be even and assume one of the values $n=2m= 2,4,8,16$. Moreover, the homotopy type of the manifold turns out to vary among a finite number of cases, including (or reducing to, if $n=2,4$) the homotopy tupe of the projective plane over the real, complex, quaternion or Cayley numbers.\\
{\bf{Definition \ref{reeb}.10}} In view of the results of Eells and Kuiper \cite{Ku-Ee}, if a manifold $M$ admits a real-valued function with exactly three non-degenerate singular points, we will say that $M$ is an {\em{Eells-Kuiper}} manifold.\\
We have (see \cite{Ca-Sca} for the case $n=3$):\\
{\bf{Theorem \ref{reeb}.11} (Center-Saddle Theorem)} {\em{Let $M^n$ be an $n$-dimensional manifold, with $n \geq 2$ such that $(M,\F)$ is a foliated manifold, by means of a transversely orientable, codimension-one, Morse, $C^\infty$ foliation $\F$. Moreover $\F$ is assumed to be without holonomy if $n=2$. Let $Sing(\F)$ be the singular set of $\F$, with $\# Sing(\F)=k+l$, where $k,l$ are the numbers of, respectively, centers and saddles. If we have $k \geq l+1$, then there are two possibilities:\\
{\em{(1)}} $k=l+2$ and $M$ is homeomorphic to an $n$-dimensional sphere;\\
{\em{(2)}} $k=l+1$ and $M$ is an Eells-Kuiper manifold}}.\\
{\em{Proof.}} If $l=0$, assertion is proved by the Reeb Sphere Theorem \ref{uno}.4. Let $l \geq 1$; we prove our thesis by induction on the number $l$ of saddles. We set $\F_l=\F$.\\
So let $l=1$ and $\F_1=\F$. By hypothesis, in the set $Sing(\F)$ there exist at least two centers, $p_1,p_2$, with $p_1 \neq p_2$, and one saddle $q$. We have necessarily $q \in \partial \C_{p_1}(\F) \cap \partial \C_{p_2}(\F)$. In fact, if this is not the case and, for example $q \notin \partial \C_{p_1}(\F)$, then (keeping into account that for $n=2$, the foliation $\F$ is assumed to be without holonomy) $\partial \C_{p_1}=\varnothing$ and $M=\overline{\C_{p_1}(\F)}$. A contraddiction. Let $i(q)$ the Morse index of the saddle $q$.\\
For $n \geq 3$ we apply the results of Theorems \ref{tre}.4 and \ref{tre}.5 to the couples $(p_1,q)$ and $(p_2,q)$. In particular, by Theorem \ref{tre}.5, {\em{(iii)}}, it follows that the saddle $q$ cannot be selfconnected. We now have the following two possibilities:\\
{\em{(a)}} $i(q)=1,n-1$ and $(p_1,q)$ or (and) $(p_2,q)$ is a trivial couple,\\
{\em{(b)}} $i(q) \neq 1,n-1$ and there are no trivial couples.\\
For $n=2$, we have necessarily $i(q)=1$ and, in our hypotheses, $q$ is always selfconnected. With few changes, we adapt Theorem \ref{tre}.5, to this case, obtaining $\partial \C_p(\F) \simeq S^1$ or $\partial \C_p(\F) \simeq S^1 \vee S^1$; in this latter case we will say that the saddle $q$ is {\em{selfconnected with respect to}} $p$. We obtain:\\
{\em{(a')}} $(p_1,q)$ or (and) $(p_2,q)$ is a trivial couple;\\
{\em{(b')}} $q$ is selfconnected both with respect to $p_1$ and to $p_2$.\\
In cases {\em{(a)}} and {\em{(a')}} we proceed with the elimination of a trivial couple, as stated in Proposition \ref{quattro}.2, and then we obtain the foliated manifold $(M,\F_0)$, with no saddle-type and some center-type singularities. We apply the Reeb Sphere Theorem \ref{uno}.4 and obtain $\# Sing(\F)=2$ and $M \simeq S^n$.\\
In case {\em{(b)}} ($n \geq 3$), as a consequence of Theorem \ref{tre}.4, we necessarily have $i(q)=n/2$ (and therefore $n$ must be even!). Moreover $\overline{\C_{p_1}(\F)} \approx \overline{\C_{p_2}(\F)}$ and $M=\overline{\C_{p_1}(\F)} \cup_\phi \overline{\C_{p_2}(\F)}$ may be thought as two copies of the same (singular) manifold glued together along the boundary, by means of the diffeomorphism $\phi$.\\
In case {\em{(b')}} ($n=2$), we obtain the same result as above, i.e. $\overline{\C_{p_1}(\F)} \approx \overline{\C_{p_2}(\F)}$ and $M=\overline{\C_{p_1}(\F)} \cup_\phi \overline{\C_{p_2}(\F)}$. We notice that case {\em{(b')}} occurs when the set $\C_{p_i}(\F) \simeq \bo ^2/S^0$ is obtained by identifying two points of the boundary in a way that reverses the orientation.\\
In cases {\em{(b)}} and {\em{(b')}}, it turns out that $\# Sing (\F_1)=3$. Moreover, $\F_1$ has a first integral, which is given by the projection of $M$ onto the space of (possibly singular) leaves. In fact,  by Lemma \ref{reeb}.6, the space of leaves is diffeomorphic to a closed interval of $\R$. In this way $M$ turns out to be an Eells-Kuiper manifold. This ends the case $l=1$.\\
Let $l>1$ (and $\# Sing(\F) >3$). As above, in $Sing(\F)$ there exist at least one saddle $q$ and two (distinct) centers, $p_1, p_2$ such that $q \in \partial \C_{p_1}(\F) \cap \partial \C_{p_2}(\F)$; we are led to the same possibilities  {\em{(a)}}, {\em{(b)}} for $n \geq 3$ and {\em{(a)'}}, {\em{(b)'}} for $n=2$. Anyway {\em{(b)}} and {\em{(b')}} cannot occur, otherwise $M=\overline{\C_{p_1}(\F)} \cup_\phi \overline{\C_{p_2}(\F)}$ and $\# Sing(\F)=3$, a contraddiction. Then we may proceed with the elimination of a trivial couple. In this way we obtain the foliated manifold $(M,\F_{l-1})$, which we apply the inductive hypothesis to. The theorem is proved, observing that, a posteriori, case {\em{(1)}} holds if $k=l+2$ and case {\em{(2)}} if $k=l+1$.\\
\section{Haefliger-type theorems}\label{haefli}
In this paragraph, we investigate the existence of leaves of singular foliations with unilateral holonomy. Keeping into account the results of the previous paragraph, for Morse foliations, we may state or exclude such an occurrence, according to the following theorem:\\
{\bf{Theorem \ref{haefli}.1}} {\em{Let $\F$ be a $C^ \infty$, codimension one, Morse foliation on
a compact manifold $M^n$, $n \geq 3$, assumed to be transversely orientable, but not necessarily closed. Let $k$ be the number of centers and $l$ the number of saddles. We have the
following possibilities: {\em{(i)}} if $k \geq l+1$, then all leaves are closed in $M \setminus Sing(\F)$; in particular, if $\partial M \neq \varnothing$ or $k \geq l+2$ each regular (singular) leaf of $\F$, is diffeomorphic (homeomorphic) to a sphere (in the second option, it is diffeomorphic to a sphere with a pinch at one point); {\em{(ii)}} if $k=l$ there are two possibilities: all leaves are closed in $M \setminus Sing(\F)$, or
there exists some compact (regular or singular) leaf with unilateral
  holonomy}}.\\
{\bf{Example \ref{haefli}.2}} The foliation of example \ref{quattro}.7 is an occurrence of theorem \ref{haefli}.1, case {\em{(ii)}} with all leaves closed. The Reeb foliation of $S^3$ and each foliation we may obtain from it, with the introduction of $l=k$ trivial couples center-saddle, are examples of theorem \ref{haefli}.1, case {\em{(ii)}}, with a leaf with unilateral holonomy.  

Now we consider other possibilities for $Sing(\F)$.\\
{\bf{Definition \ref{haefli}.3}} Let $\F$ be a $C^ \infty$, codimension one foliation on a compact
manifold $M^n, n\geq 3$, with singular set $Sing(\F) \neq \varnothing$. We
say that $Sing(\F)$ is {\em{regular}} if its connected components
are either isolated points or smoothly embedded curves,
diffeomorphic to $S^1$. We extend the definition of stability to regular components, by saying that a connected component $\Gamma \subset Sing(\F)$ is {\em{(weakly) stable}}, if there exists a neighborhood of $\Gamma$, where the foliation has all leaves compact (notice that we can repeat the proof of Lemma \ref{reeb}.7 and obtain that a weakly stable component is a stable component).

In the case $Sing(\F)$ is regular, with stable isolated singularities, when $n \geq 3$ we may exclude a Haefliger-type result, as a consequence of Lemma \ref{reeb}.5 and the Reeb Global Stability Theorem for manifolds with boundary. Then we study the case $Sing(\F)$ regular, with stable components, all diffeomorphic to $S^1$. Let $J$ be a set such that for all $j \in J$, the curve $\gamma_j:S^1 \rightarrow M$, is a smooth embedding and $\Gamma_j:= \gamma_j(S^1) \subset Sing(\F)$ is stable. Then $J$ is a finite set. This is obvious, otherwise $\forall j \in J$, we may select a point $x_j \in \Gamma_j$ and obtain that the set $\{x_j \}_{j \in J}$ has an accumulation point. But this is not possible because the singular components are separated. We may regard a singular component $\Gamma_j$, as a {\em{degenerate leaf}}, in the sense that we may associate to it, a single point of the space of leaves.

We need the following definition\\
{\bf{Definition \ref{haefli}.4}} Let $\F$ be a $C^ \infty$, codimension one foliation on a compact manifold $M$. Let $\overline{D^2}$ be the closed 2-disc and $g:\overline{D^2} \rightarrow M$ be a $C^ \infty$ map. We say that $p \in \overline{D^2}$ is a {\em{tangency point of $g$ with $\F$}} if
$(\textrm{d}g)_p (\R^2) \subset T_{g(p)} \F_{g(p)}.$

We recall a proposition which Haefliger's theorem (cfr. the book \cite{Cam}) is based upon.\\ 
{\bf{Proposition \ref{haefli}.5}} {\em{Let $A: \overline{D^2} \rightarrow M$ be a $C^ \infty$ map, such that the restriction $A|_{\partial D^2}$ is transverse to $\F$, i.e. $\forall x \in \partial D^2,(\textrm{d}A)_x(T_x (\partial D^2))+ T_{A(x)} \F_{A(x)}=T_{A(x)}M$. Then, for every $\epsilon >0$ and every integer $r \geq 2$, there exists a $C^ \infty$ map, $g: \overline{D^2} \rightarrow M$, $\epsilon$-near $A$ in the $C^r$-topology, satisfying the following properties: {\em{(i)}} $g|_{\partial D^2}$ is transverse to $\F$. {\em{(ii)}}
For every point $p \in D^2$ of tangency of $g$ with $\F$, there
exists a foliation box $U$ of $\F$ with $g(p) \in U$ and a
distinguished map $\pi:U \rightarrow \R$ such that $p$ is a
non-degenerate singularity of $\pi \circ g:g^{-1}(U) \rightarrow
\R$. In particular there are only a finite number of tangency
points of $g$ with $\F$, since they are isolated, and they are contained in the
open disc $D^2=\{z \in \R^2:||z||<1\}$. {\em{(iii)}} If $T=\{p_1,
\dots, p_t \}$ is the set of tangency points of $g$ with $\F$,
then $g(p_i)$ and $g(p_j)$ are contained in distinct leaves of
$\F$, for every $i \neq j$. In particular, the singular foliation
$g^*(\F)$ has no saddle connections.}}

We are now able to prove a similar result, in the case of existence of singular components.\\{\bf{Proposition \ref{haefli}.6}} {\em{Let $\F$ be a codimension one, $C^ \infty$ foliation on a compact manifold $M^n$, $n \geq 3$, with regular singular set, $Sing(\F) = \cup _{j \in J} \Gamma_j \neq \varnothing$, where $\Gamma_j$ are all stable components diffeomorphic to $S^1$ and $J$ is finite. Let $A: \overline{D^2} \rightarrow M$ be a $C^ \infty$ map, such that the restriction $A|_{\partial D^2}$ is transverse to $\F$. Then, for every $\epsilon >0$ and every integer $r \geq 2$, there exists a $C^ \infty$ map, $g: \overline{D^2} \rightarrow M$, $\epsilon$-near $A$ in the $C^r$-topology, satisfying properties {\em{(i)}} and {\em{(iii)}} of proposition \ref{haefli}.5, while {\em{(ii)}} is changed in: {\em{(ii')}} for every point  $p \in D^2$ of tangency of $g$ with
$\F$, we have two cases: (1) if $L_{g(p)}$ is a regular leaf of $\F$, there exists a foliation box, $U$ of $\F$, with $g(p) \in U$, and a distinguished map, $\pi:U \rightarrow \R$, satisfying properties as in {\em{(ii)}} of Proposition \ref{haefli}.5; (2) if $L_{g(p)}$ is a degenerate leaf of $\F$, there exists a neighborhood, $U$ of $p$, and a singular submersion, $\pi:U \rightarrow \R$, satisfying properties as in {\em{(ii)}} Proposition \ref{haefli}.5.}}\\
{\em{Proof.}} We start by recalling the idea of the classical proof.\\
We choose a finite covering of $A(\overline{D^2})$ by foliation boxes $\{Q_i\}^r_{i=1}$. In each $Q_i$ the foliation is defined by a distinguished map, the submersion $\pi_i:Q_i \rightarrow \R$. We choose an atlas, $\{(Q_i, \phi_i)\}
^r_{i=1}$, such that the last component of $\phi_i:Q_i \rightarrow \R^n$ is $\pi_i$, i.e. $\phi_i=(\phi_i^1, \phi_i^2, \dots , \phi_i^{n-1}, \pi_i)$. We construct the finite cover of $\overline{D^2}$, $\{W_i=A^{-1}(Q_i)\}^r_{i=1}$; the expression of $A$ in coordinates is $A|_{W_i}=(A_i^1, \dots , A_i^{n-1}, \pi_i \circ A)$. We may choose covers of $\overline{D^2}$, $\{U_i \}_{i=1}^r$, $\{V_i \}_{i=1}^r$, such that  $\overline{U_i} \subset V_i \subset \overline {V_i} \subset W_i$, $i=1, \dots , r$; then we proceed by induction on the number $i$. Starting with $i=1$ and setting $g_0=A$, we apply a result (\cite{Cam}, Cap. VI, $\S$2, Lemma 1, pag. 120) and we modify $g_{i-1}$ in a new function $g_i$, in a way that $g_i(W_i) \subset Q_i$ and $\pi_i \circ g_i:W_i \rightarrow \R$ is Morse on the subset $U_i \subset W_i$. At last we set $g=g_r$.

In the present case, essentially, it is enough to choose a set of couples, $\{(U_k,\pi_k)\}_{k \in K}$, where $\{U_k\}_{k \in K}$ is an open covering of $M$, $\pi_k:U_k \rightarrow \R$, for $k \in K$, is a (possibly singular) submersion and, if $U_k \cap U_l \neq \varnothing$ for a couple of indices $k,l \in K$, there exists a diffeomorphism $p_{lk}:\pi_k(U_k \cap U_l) \rightarrow \pi_l(U_k \cap U_l)$, such that $\pi_l=p_{lk}\circ \pi_k$. By hypothesis, there exists the set of couples $\{(U_i,\pi_i)\}_{i \in I}$, where $\{U_i\}_{i \in I}$, is an open covering of $M \setminus Sing(\F)$, and, for $i \in I$, the map $\pi_i:U_i \rightarrow  \R$, is a distinguished map, defining the foliated manifold $(M \setminus Sing(\F), \F^*)$. Let $y \in Sing(\F)$, then $y \in \Gamma_j$, for some $j \in J$. As $y \in M$, there exists a neighborhood $C \ni y$, homeomorphic to an $n$-ball. Let $h:C \rightarrow \bo^n$ be such a homeomorphism. As the map $\gamma_j:S^1 \rightarrow \Gamma_j$ is a smooth embedding, we may suppose that, locally, $\Gamma_j$ is sent in a diameter of the ball $\bo^n$, i.e. $h(C \cap \Gamma_j)=\{x_2=\dots =x_n=0 \}$. For each singular point $z=h^{-1}(b,0,\dots,0)$, the set $D=h^{-1}(b,x_2,\dots ,x_n)$, homeomorphic to a small $(n-1)$-ball, is transverse to the foliation at $z$. Moreover, if $z_1 \neq z_2$, then $D_1 \cap D_2=\varnothing$. The restriction $\F|_D$ is a singular foliation with an isolated stable singularity at $z$. By lemma \ref{reeb}.5, the leaves of $\F|_D$ are diffeomorphic to $(n-2)$-spheres. It turns out that $y$ has a neighborhood homeomorphic to the product $(-1,1) \times \bo^{n-1}$, where the foliation is the image of the singular trivial foliation of $(-1,1) \times \bo^{n-1}$, given by $(-1,1) \times S^{n-2}\times \{t \}, t \in (0,1)$, with singular set $(-1,1) \times \{0 \}$. Let $\pi_y:U_y \rightarrow[0,1)$ be the projection. If, for a couple of singular points $y,w \in Sing(\F)$, we have $U_y \cap U_w \neq \varnothing$, we may suppose they belong to the same connected component, $\Gamma_j$. We have $\pi_w \circ \pi_y^{-1}(0)=0$ and, as a consequence of lemma \ref{reeb}.6, there exists a diffeomorphism between $\pi_y(U_y \cap U_w \setminus \Gamma_j)$ and $\pi_w(U_y \cap U_w \setminus \Gamma_j)$. The same happens if $U_y \cap U_i \neq \varnothing$ for some $U_i \subset M \setminus Sing(\F)$. It comes out that $\pi_y$ is singular on $U_y \cap Sing(\F)$ and non-singular on $U_y \setminus Sing(\F)$, i.e. $(d \pi_y)_z=0 \Leftrightarrow z \in U_y \cap Sing(\F)$. At the end, we set $K=I \cup Sing(\F)$.\\
Let $g: \overline{D^2} \rightarrow M$ be a map. Then $g$ defines the foliation $g^*(\F)$, pull-back of $\F$, on $\overline{D^2}$. Observe that if $Sing(\F)= \varnothing$, then $Sing(g^*(\F))=\{ \textrm{tangency points of }g \textrm{ with } \F \}$, but in the present case, as $Sing(\F) \neq \varnothing$, we have $Sing(g^*(\F))=\{ \textrm{tangency points of }g \textrm{ with } \F \} \cup g^*(Sing(\F))$. Either if $p$ is a point of tangency of $g$ with $\F$ or if $p \in g^*(Sing(\F))$, we have $d(\pi_k)_p=0$. With this remark, we may follow the classical proof.

As a consequence of proposition \ref{haefli}.6, we have:\\
{\bf{Theorem \ref{haefli}.7 (Haefliger's theorem for singular foliations)}} {\em{Let $\F$ be a codimension one, $C^2$, possibly singular foliation of an $n$-manifold $M$, with $Sing(\F)$, (empty or) regular and with stable components diffeomorphic to $S^1$. Suppose there exists a closed curve transverse to $\F$, homotopic to a point. Then there exists a leaf with unilateral holonomy.}}
\section{Novikov-type theorems}\label{novi}
We end this article with a result based on the original Novikov's
Compact Leaf Theorem and on the notion of stable
singular set. To this aim, we premise the following remark. Novikov's statement establishes the existence of a compact leaf for foliations on 3-manifolds with finite fundamental group. This result actually proves the existence of an invariant submanifold, say $N \subset M$, with boundary, such that $\F|_N$ contains open leaves whose universal covering is the plane. Moreover these leaves accumulate to the compact leaf of the boundary. In what follows, a submanifold with the above properties will be called a {\em{Novikov component}}. In particular a Novikov component may be a Reeb component, i.e. a solid torus endowed with its Reeb foliation. We recall that two Reeb components, glued together along the boundary by means of a diffeomorphism which sends meridians in parallels and viceversa, give the classical example of the Reeb foliation of $S^3$.\\
If $\F$ is a Morse foliation of a 3-manifold, as all saddles have index 1 or 2, we are always in conditions of proposition \ref{quattro}.2 and then we are reduced to consider just two (opposite) cases: {\em{(i)}} all singularities are centers, {\em{(ii)}} all singularities are saddles. In case {\em{(i)}}, by the proof of the Reeb Sphere Theorem \ref{uno}.4, we know that all leaves are compact; in case {\em{(ii)}}, all leaves may be open and dense, as it is shown by an example of a foliation of $S^3$ with Morse singularities and no compact leaves \cite{Ros-Rou}.\\
As in the previous paragraph, we study the case in which $Sing(\F)$ is regular with stable components, $\Gamma_j, j \in J$, where $J$ is a finite set. We have:\\
{\bf{Theorem \ref{novi}.1}} {\em{Let $\F$ be a $C^ \infty$, codimension one foliation on a closed
$3$-manifold $M^3$. Suppose
$Sing(\F)$ is (empty or) regular, with stable components. Then we have
two possibilities: {\em{(i)}} all leaves of $\F$ are compact; {\em{(ii)}} $\F$ has a Novikov
component.}}\\
{\em{Proof.}} If $Sing(\F)= \varnothing$, thesis (case {\em{(ii)}}) follows by Novikov theorem. Let $Sing(\F) \neq \varnothing$. We may suppose that $\F$ is transversely orientable (otherwise we pass to the transversely orientable double covering). If $Sing(\F)$ contains an isolated singularity, as we know, we are in case {\em{(i)}}. Then we suppose $Sing(\F)$ contains no isolated singularity, i.e. $Sing(\F)= \bigcup_{j \in J} \Gamma_j$. Set $\mathcal D(\F)= \{\Gamma_j, j \in J \} \cup \{ \textrm{ compact leaves with trivial holonomy} \}$. By the Reeb Local Stability Theorem \ref{uno}.1,  $\mathcal D(\F)$ is open. We may have $\partial \mathcal D(\F)= \varnothing$, and then we are in case {\em{(i)}}, or $\partial \mathcal D(\F)\neq \varnothing$, and in this case it contains a leaf with unilateral holonomy, $F$. It is clear that $F$ bounds a Novikov component, and then we are in case {\em{(ii)}}; in fact, from one side, $F$ is accumulated by open leaves. If $F'$ is one accumulating leaf, then its universal covering is $p:\R^2 \rightarrow F'$. Suppose, by contraddiction, that the universal covering of $F'$ is $p:S^2 \rightarrow F'$. By the Reeb Global Stability Theorem for manifolds with boundary, all leaves are compact, diffeomorphic to $p(S^2)$. This concludes the proof since $F$ must have infinite fundamental group.

The last result may be reread in terms of the existence of closed curves, transverse to the foliation. We have:\\
{\bf{Lemma \ref{novi}.2}} {\em{Let $\F$ be a codimension one, $C^ \infty$ foliation on a
closed $3$-manifold $M$, with singular set, $Sing(\F) \neq \varnothing$, regular, with 
stable components. Then $\F$ is a foliation with all leaves compact if and only if there exist no closed transversals.}}\\
{\em{Proof.}} (Sufficiency) If the foliation admits an open (in $M \setminus Sing(\F)$) leaf, $L$, it is well known that we may find a closed curve, intersecting $L$, transverse to the foliation. Viceversa (necessity), let $\F$ be a foliation with all leaves compact. If necessary, we pass to the transversely orientable double covering $p:(\widetilde{M}, \widetilde{\F}) \rightarrow (M,\F)$. In this way, we apply Lemma \ref{reeb}.6 and obtain, as $Sing(\widetilde{\F}) \neq \varnothing$, that the projection onto the space of leaves is a (global) $C^ \infty$ first integral of $\widetilde{\F}$, $f:\widetilde{M} \rightarrow [0,1] \subset \R$. Suppose, by contraddiction, that there exists a $C^1$ closed transversal to the foliation $\F$, the curve $\gamma:S^1 \rightarrow M$. The lifting of $\gamma^2$ is a closed curve, $\Gamma:S^1 \rightarrow \widetilde{M}$, transverse to $\widetilde{\F}$. The set $f(\Gamma(S^1))$ is compact and then has maximum and minimum, $m_1,m_2 \in \R$. A contraddiction, because $\Gamma$ cannot be transverse to the leaves $\{ f^{-1}(m_1) \}, \{ f^{-1}(m_2)\}$.

\begin{figure}
\begin{minipage}[t]{.45\linewidth}
  {
\begin{center}
{
        \psfrag{p}{$p$}
        \psfrag{q}{$q$}
        \psfrag{L}{$L$}
        \psfrag{F}{$F$}
        \psfrag{U}{$U$}
        \psfrag{S}{$\partial \C_p(\F)$}
        \includegraphics[scale=.25]{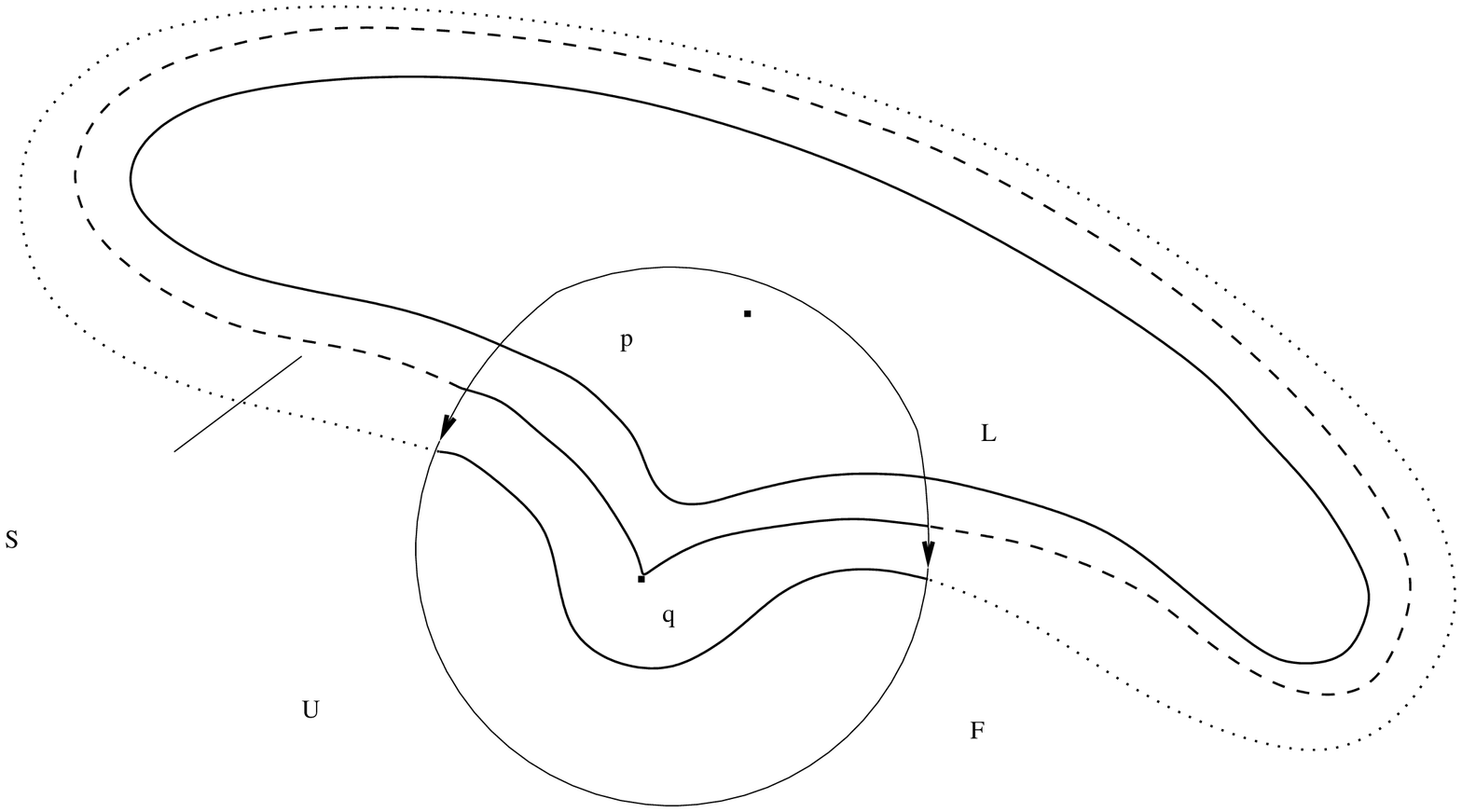}
        \caption{$p-q$ is not a trivial coupling when $1<l<n-1$, where $l$ is
        the index of the saddle $q$.} 
        \label{LF}    
}
\end{center}
  }\end{minipage}%
  \begin{minipage}[t]{.1\linewidth}{\hspace{.1\linewidth}}\end{minipage}%
  \begin{minipage}[t]{.45\linewidth}
{
\begin{center}
{
        \psfrag{1}{$ST_1$}
        \psfrag{2}{$ST_2$}
        \psfrag{g}{$\gamma$}
        \includegraphics[scale=.35]{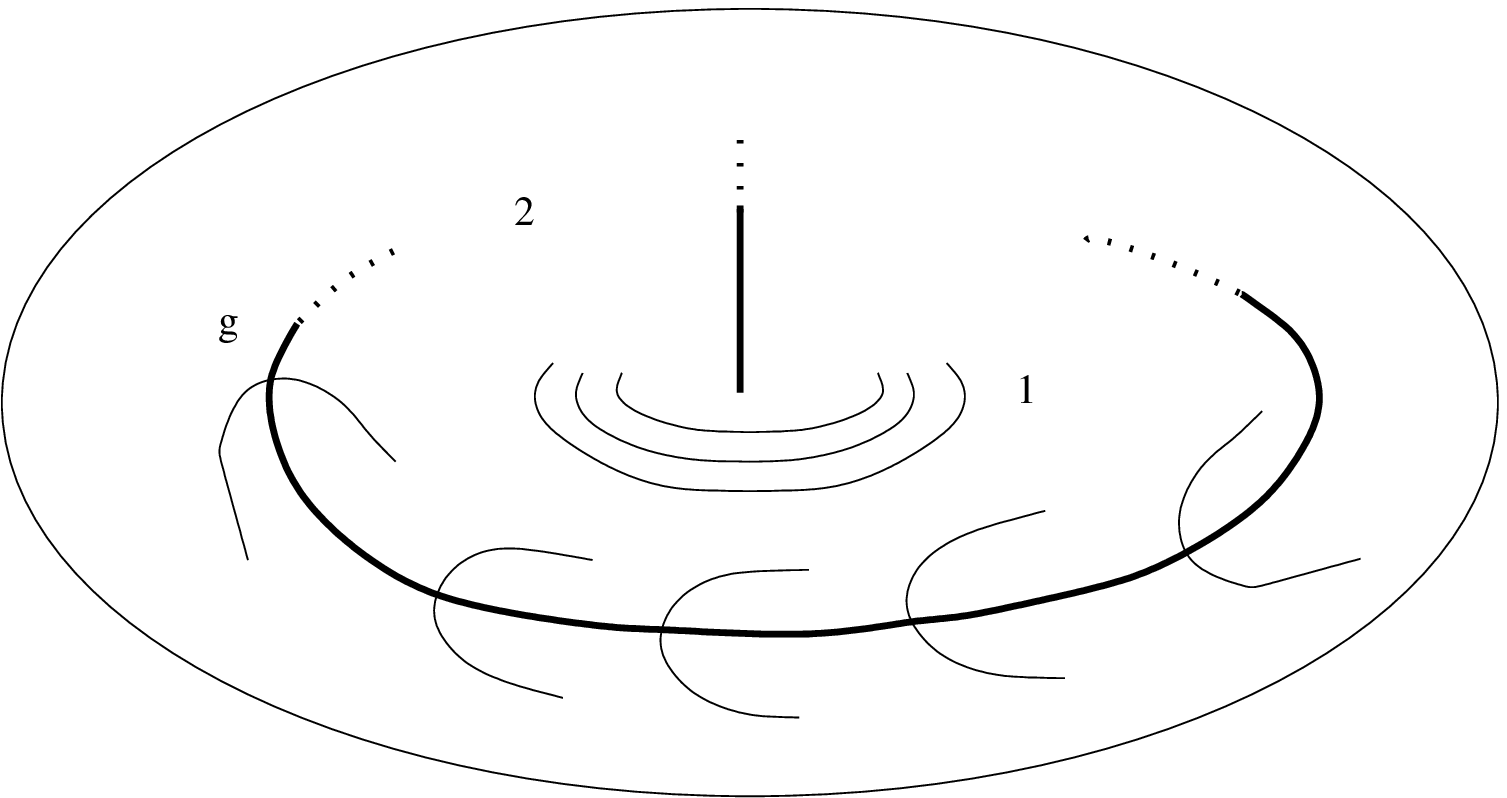}
        \caption{A singular foliation of $S^3$, with no vanishing cycles.}
 \label{nocycle}
}
\end{center}
}
\end{minipage}
\end{figure}
With this result, we may rephrase the previous theorem.\\
{\bf{Corollary \ref{novi}.3}} {\em{Let $\F$ be a codimension one, $C^ \infty$ foliation on a $3$-manifold $M$, such that $Sing(\F)$ is regular with stable components. Then {\em{(i)}} there are no closed transversals, or equivalently, $\F$ is a foliation by compact leaves, {\em{(ii)}} there exists a closed transversal, or equivalently, $\F$ has a Novikov component.}}\\
{\bf{Remark \ref{novi}.4}} In the situation we are considering, we cannot state a singular version of Auxiliary Theorem I (see, for example \cite{Mo-Sca}). In fact, even though a singular version of Haefliger Theorem is given, the existence of a closed curve transverse the foliation, homotopic to a constant, does not lead, in general, to the existence of a vanishing cycle, as it is shown by the following counterexample.\\
{\bf{Example \ref{novi}.5}} We consider the foliation of $S^3$ given by a Reeb component, $ST_1$, glued (through a diffeomorphism of the boundary which interchanges meridians with parallels) to a solid torus $ST_2= S^1 \times \overline{D^2}=T^2 \times (0,1) \cup S^1$. The torus $ST_2$ is endowed with the singular trivial foliation $\F|_{ST_2}=T^2 \times \{t \}$, for $t \in (0,1)$, where $Sing(\F|_{ST_2})=S^1=Sing(\F)$. As a closed transversal to the foliation, we consider the curve $\gamma:S^1 \rightarrow ST_1 \subset S^3$, drawed in figure \ref{nocycle}. Let $f:\overline{D^2} \rightarrow S^3$ be an extension of $\gamma$; the extension $f$ is assumed to be in general position with respect to $\F$, as a consequence of proposition \ref{haefli}.5. As $\gamma(S^1)$ is linked to the singular component $S^1 \subset ST_2$, then $f(\overline{D^2}) \cap Sing(\F) \neq \varnothing$. As a consequence, we find a decreasing sequence of cycles, $\{\beta_n \}$, (the closed curves of the picture) which does not admit a cycle, $\beta_\infty$, such that $\beta_n > \beta _\infty$, for all $n$. In fact the ``limit'' of the sequence is not a cycle, but the point $f(\overline{D^2}) \cap Sing(\F)
$. \\
{\bf{Example \ref{novi}.6}} The different situations of Theorem \ref{novi}.1 or Corollary \ref{novi}.3 may be exemplified as follows. It is easy to see that $S^3$ admits a singular foliation with all leaves compact (diffeomorphic to $T^2$) and two singular (stable) components linked together, diffeomorphic to $S^1$. In fact one can verify that $S^3$ is the union of two solid tori, $ST_1$ and $ST_2$, glued together along the boundary, both endowed with a singular trivial foliation.\\
We construct another foliation on $S^3$, modifying the previous one. We set $\widetilde{ST_1}=S^1 \times \{0 \} \cup T^2 \times (0,1/2]$. In this way, $ST_1= \widetilde{ST_1} \cup T^2 \times (1/2,1]$. We now modify the foliation in $ST_1 \setminus \widetilde{ST_1}$, by replacing the trivial foliation with a foliation with cylindric leaves accumulating to the two components of the boundary.


\begin{thebibliography}{99}
\bibitem[Cam-LN]{Cam} C. Camacho, A. Lins Neto: Geometric theory of foliations, Boston, Birkhauser, 1985
\bibitem[Cam-Sc]{Ca-Sca} C. Camacho, B. Sc\'ardua: On codimension one foliations with Morse singularities on three-manifolds, Topology and its Applications 154 (2007) 1032-1040.
\bibitem[Ee-Kui]{Ku-Ee} J. Eells, N.H. Kuiper: Manifolds which are
  like projective planes, Pub. Math. de l'I.H.E.S., 14, 1962.
\bibitem[God]{God} C. Godbillon: Feuilletages, etudies geometriques, Basel, Birkhauser, 1991
\bibitem[Law]{Law} H.B. Lawson, jr.: Foliations, Bull. Amer. Math. Soc., Vol. 80, N. 3, May 1974.
\bibitem[Mil 1]{Mil1} J. Milnor: Morse theory, Princeton, NJ, Princeton
University Press, 1963.
\bibitem[Mil 2]{Mil2} J. Milnor: Lectures on the h-cobordism theorem, Princeton, NJ, Princeton University Press, 1965.
\bibitem[Mor-Sc]{Mo-Sca} C.A. Morales, B. Sc\'ardua: Geometry and Topology of foliated manifolds.
\bibitem[Nov]{Nov} S.P. Novikov: Topology of foliations. Trudy Moskov. Mat. Obshch. 14 (1965), 248-278.
\bibitem[Pal-deM]{Palis} J. Palis, jr., W. de Melo: Geometric theory of dinamical systems: an introduction, New-York, Springer,1982.
\bibitem[Reeb]{Reeb} G. Reeb: Sur les points singuliers d'une forme de
Pfaff compl\`etement int\'egrable ou d'une fonction num\'erique. CRAS 222 (1946), 847-849.
\bibitem[Ros-Rou]{Ros-Rou} H. Rosemberg, R. Roussarie: Some remarks on stability of foliations, J. Diff. Geom. 10, 1975, 207-219.
\bibitem[Stee]{Stee} N. Steenrod: The topology of fiber
  bundles, Princeton, NJ, Princeton University Press, 1951
\bibitem[Thu]{Thu} W.P. Thurston: Three-dimensional geometry and
  topology, Princeton, NJ, Princeton University Press, 1997.
\bibitem[Wag]{Wag} E. Wagneur: Formes de Pfaff \`a singularit\'es
  non d\'eg\'en\'er\'ees, Annales de l'institut Fourier, tome 28 n. 3
  (1978), p. 165-176.
\end{thebibliography}
\end{document}